\documentclass[12pt]{amsart}
\usepackage[centertags]{amsmath}

\usepackage{amsfonts}

\usepackage{amssymb}

\title{Cofinal Extensions and Coded Sets}
\author{James H. Schmerl}

\date{\today}
\def\into{\longrightarrow}

\def\pa{{\sf PA}}

\def\ISig{\rm{I}\Sigma}
\def\BSig{\rm{B}\Sigma}
\def\IDel{\rm{I}\Delta}

\def\rca{{\sf RCA}_0}

\def\wkl{{\sf WKL}_0}
\def\aca{{\sf ACA}_0}

\def\wkls{\wkl^*}

\def\rcas{\rca^{*}}

\newcommand{\half}{\frac{1}{2}}

\newcommand{\KK}{{\mathcal K}}

\newcommand{\MM}{{\mathcal M}}

\newcommand{\NN}{{\mathcal N}}

  \DeclareMathOperator{\ssy}{SSy}

        \DeclareMathOperator{\dcf}{dcf}
         \DeclareMathOperator{\gcis}{GCIS}
           \DeclareMathOperator{\GCIS}{GCIS}
   
    \DeclareMathOperator{\dom}{dom}

\def\lhdeq{\trianglelefteq}
\def\rhdeq{\trianglerighteq}
\def\rhd{\vartriangleright}
\def\qqed{\vspace{-17pt}\begin{flushright}$\square \, \square$\end{flushright} \vspace{17pt}}

 \DeclareMathOperator{\cf}{cf}
 \DeclareMathOperator{\cod}{Cod}
  \DeclareMathOperator{\Cod}{\cod}

  \DeclareMathOperator{\Th}{Th} 
  
  \DeclareMathOperator{\Def}{Def}

  \DeclareMathOperator{\Lt}{Lt}

\begin{document}
\maketitle 

\begin{abstract}Let $\MM$ be a model of Peano Arithmetic that is countably generated over an exponentially closed cut $I$. We will characterize (in Theorems~2.1 and~3.3) those ${\mathfrak X} \subseteq \mathcal{P}(I)$ such that there is a finitely (or countably) generated extension $\NN \succ_{\sf cf} \MM$ such that $\gcis(\MM,\NN) = I$ and $\cod(\NN /I) = {\mathfrak X}$. We also characterize such ${\mathfrak X}$ for which the extension $\NN$ can be, in addition, one of the following: non-filling (Theorem~5.6); filling (Theorem~4.1); $n$-filling (Corollary~4.5). \end{abstract}

\bigskip

The fundamental MacDowell-Specker Theorem asserts that every model $\MM$ of Peano Arithmetic has 
an elementary end extension $\NN \succ_{\sf end} \MM$ (which clearly can be a finitely generated extension). Phillips \cite{ph} observed that this theorem could be improved by requiring $\NN$ to be a conservative extension. For any  $\NN \models \pa$, the set of its parametrically definable sets is 
$\Def(\NN)$, and if $X \subseteq N$, then the set of coded subsets of $X$ is 
$$
\cod(\NN / X) = \{X \cap D : D \in \Def(\NN)\}.
$$ 
The extension 
$\NN \succ \MM$ is {\bf conservative} if  $\cod(\NN / M) = \Def(\MM)$. Gaifman \cite{g} improved the MacDowell-Specker Theorem in another direction by proving that 
every model $\MM \models \pa$  has 
a minimal elementary end extension $\NN$, where the extension $\NN \succ \MM$ is {\bf minimal} if there is no $\KK$ such that $\MM \prec \KK \prec \NN$.  Phillips \cite{ph} also noted  that Gaifman's theorem could be improved by requiring $\NN$ to be a conservative extension. We have little to say in this paper  about minimal extensions.

The MacDowell-Specker-Phillips Theorem was generalized in \cite{code} by characterizing those ${\mathfrak X} \subseteq {\mathcal P}(M)$ for which $\MM$ has a finitely (or countably) generated, elementary end extension $\NN \succ_{\sf end} \MM$ such that $\cod(\NN / M) = {\mathfrak X}$. This generalization appears  here as Corollary~2.8. 
Roman Kossak \cite{Q} proposed  the possibility that this theorem in \cite{code} could be generalized to $\aleph_1$-generated elementary end extensions. This proposal is proved here in Corollary~2.10. The proof relies on Theorem~2.1 which, together with Theorem~ 3.3, characterizes (under some additional hypotheses) the possible sets of sets that can be coded in elementary cofinal extensions. The proof of Kossak's proposal 
is the impetus for the focus of this paper, which  is to investigate which sets of sets can be coded in cofinal extensions of models of \pa. 

Given $\MM \prec \NN \models \pa$, we let  their 
\underline{G}reatest \underline{C}ommon \underline{I}nitial \underline{S}egment  be the set 
$$
\gcis(\MM, \NN) = \{b \in M : [0,b]^\MM = [0,b]^\NN\},
$$
 where $[a,b]^\MM = \{x \in M : \MM \models a \leq x \leq b\}$. If $\MM \prec_{\sf cf} \NN$, then 
 $\gcis(\MM, \NN)$ is a cut of both $\MM$ and $\NN$. 
 Every cut $\gcis(\MM, \NN)$ is closed under multiplication (Proposition 1.2(1)).  For countable $\MM$, there is a converse \cite[Th.~2]{pm}: if $I$ is a cut of $\MM$ closed under multiplication, then there is $\NN \succ_{\sf cf} \MM$ such that 
 $\gcis(\MM, \NN) = I$.  Nevertheless, this paper will usually consider only  cuts closed under exponentiation.\footnote{I leave the study of coded sets in extensions $\NN \succ_{\sf cf} \MM$, where $\gcis(\MM, \NN)$ may not be closed under exponentiation, to a later time or another researcher.} If $\MM \prec_{\sf cf} \NN$, then we will be concerned with the set $\cod(\MM / \gcis(\MM, \NN))$, which is what we think of as the set of sets coded in the cofinal extension. 
There had been some much earlier work on $\gcis(\MM, \NN)$ (for example, in \cite{pm}) in which $\MM$ was usually restricted to being countable.  Indeed, we will usually require a countability condition, but a weaker one. If $\MM \models \pa$,  
then $X$ {\bf generates} $\MM$ if  $X \subseteq M$ and there is no $\KK \prec \MM$ such 
that $X \subseteq K$.  If $X,I \subseteq M$, then $X$ {\bf generates} $\MM$ over $I$ if $X \cup I$ generates $\MM$. For any cardinal $\kappa$, we say that $\MM$ is $\kappa$-{\bf generated} over $I$ if  there is $X$ that generates $\MM$ over $I$ such that $|X| \leq \kappa$. An extension $\NN \succ \MM$ is $\kappa$-{\bf generated} if then model $\NN$ is $\kappa$-generated over $M$. We say that $\MM$ is {\bf countably} (or {\bf finitely}) generated over $I$ if it is ${\aleph_0}$-generated (or $1$-generated)  over $I$. For any $\MM \models \pa$ that  is countably generated over an exponentially closed cut $I$, we characterize (in Theorems~2.1 and~3.3)  those  ${\mathfrak X} \subseteq {\mathcal P}(I)$ for which there is 
 a finitely (or countably) generated extension $\NN \succ_{\sf cf} \MM$ such that $I = \gcis(\MM,\NN)$ and ${\mathfrak X} = 
\cod(\NN / I)$. 

Two types of cofinal extensions  will be studied here. 
 If $\MM \prec_{\sf cf} \NN \models \pa$, then $\NN$ is a {\bf filling} extension if there is $a \in N$ such that $x < a < y$ whenever $x \in \gcis(\MM,\NN)$ and $y \in M \backslash \gcis(\MM, \NN)$; otherwise, the extension is {\bf non-filling}. These types of extensions originated with Mills\cite{mills} and Paris \& Mills \cite{pm}. Theorem~5.6 deals with non-filling extensions and the possible sets of sets that can be coded in such extensions. Theorem~4.1 is its analog for  filling extensions.


\bigskip 

{\bf \S1.\@ Preliminaries.} The theory $\pa$ is formulated in the  first-order language ${\mathcal L}_\pa = \{+,\times, 0,1,\leq\}$. Script letters such as $\MM, \NN, \MM_0$, etc., will be used for ${\mathcal L}_\pa$-structures, with universes $M,N,M_0$, etc., and will always be models of some weak fragment of $\pa$ (say $\ISig_0$), but often they will be models  of \pa.  We assume that every $\MM$ is an end extension of the standard model 
${\mathbb N} = (\omega,+,\times,0,1,\leq)$. If $\NN \subseteq \MM$ (that is, $\NN$ is a submodel of $\MM)$, then $\MM | N = \NN$. 

Suppose that $\MM$ is a model that is currently under consideration. If $C \subseteq M$, then ${\mathcal L}(C)$ is  the language ${\mathcal L}_\pa$ augmented by constants denoting elements of $C$, and we let 
$T_C $ be the ${\mathcal L}(C)$-theory $\Th((\MM,c)_{c \in C})$. In particular, $T_M$ is the complete diagram of $\MM$. We let each of $\Sigma_0, \Delta_0, \Pi_0$ be the set of bounded ${\mathcal L}(M)$-formulas, and then define $\Sigma_n, \Pi_n$ ($n < \omega$) accordingly. 

 As already mentioned in the introduction, $\Def(\MM)$ is the set of those subsets of $M$ that are parametrically definable in $\MM$, and if $A \subseteq M$, then $\cod(\MM / A) = \{A \cap D : D \in \Def(\MM)\}$ is the set of coded subsets of $A$. If $\MM$ is nonstandard, then the {\bf standard system} of $\MM$ is $\ssy(\MM) = \cod(\MM / \omega)$.   If $A \subseteq M$, then $\MM|A$ is the substructure of $\MM$ having universe $A$. It is to be tacitly understood, when considering
such a substructure, that $\omega \subseteq A$ and that $A$ is closed under $+$ and $\times$. 
A {\bf cut} $I$ of $\MM$ is a proper subset of $M$ (i.e., $\varnothing \varsubsetneq I \varsubsetneq M$) 
such that $a \in I$ whenever  $b \in I$ and $a \leq b+1$. A cut $I$ is {\bf closed under} addition, multiplication, exponentiation (or is additively, multiplicatively or exponentially closed)  if, respectively,  $2a, a^2, 2^a \in I$ whenever $a \in I$. It is {\bf elementary} if $\MM | I \prec \MM$.    
If $X \subseteq M$, then $\inf(X) = \{ y \in M : y < x $ for all $x \in X\}$ and $\sup(x) = \{y \in M : y \leq x$ for all $x \in X\}$. If $I$ is a cut (or $I = M$),
then its {\bf cofinality} $\cf(I)$ is the smallest cardinal $\kappa$ for which there is $X \subseteq I$ such that $|X| = \kappa$ and $I = \sup(X)$.  Its {\bf downward cofinality} $\dcf(I)$ is the smallest cardinal $\kappa$ for which there  $X \subseteq M \backslash I$ such that $|X| = \kappa$ and $I = \inf(X)$.

In the customary way, define the ordered pair $\langle x,y \rangle$ 
to be $(x^2 + y^2 + 2xy + 3x + y)/2$, and  let exponentiation $x \mapsto 2^x$ be given  a $\Delta_0$ definition 
\`{a} la Bennett \cite{ben} (or see \cite[Chap.\@ V.3]{hp}). If $\MM$ is a model, $A \subseteq M$ and $a \in M$, then 
$(A)_a = \{x \in M : \langle a,x \rangle \in A\}$. It may be that exponentiation is a partial function; we let {\sf exp} be the statement asserting that it is total.

If $\MM \models \pa$ and $m \in M$, then an $\MM$-{\bf sequence} of {\bf length} $m$ is an $\MM$-definable function $f : [0,m-1]^\MM \into M$. We arrange that every 
$a$ in $M$ encodes a unique sequence: we let $\ell(a)$ be the length of the sequence that $a$ encodes; 
if $i < \ell(a)$, then $(a)_i$ is the $i$-th element of the sequence that $a$ encodes. If $a,b \in M$, $a \lhdeq b$ iff the sequence encoded by $a$ is an initial segment of the one encoded by $b$. These definitions are made so that they are $\Delta_0$ and  whenever   $I$ is an exponentially closed cut and $\ell(a) \in I$ and $(a)_i \in I$ for all $i < \ell(a)$, then $a \in I$. Furthermore, if $a \lhd b$, then $a < b$. 
Thus, $0$ codes the sequence of length $0$.

If $I$ is an exponentially closed cut of $\MM$ or $I = M$, then  let $2^{<I}$ be the set those $a \in M$ that code $0,1$-sequences. Thus, $2^{<I} \subseteq I$. 
 An $I$-{\bf tree} is a nonempty subset $T \subseteq 2^{<I}$ such that whenever $s \lhd t \in T$, then 
$s \in T$. If $T$ is an $I$-tree, then $P$ is a {\bf path} through $T$ or  a $T$-{\bf path}, if $P \subseteq T$, $P$ is linearly ordered by $\lhd$, and for each $t \in T$ there is $p \in P$ such that $\ell(p) = \ell(t)$. 
If $A \subseteq I$ is such that $[0,b]^\MM \cap A \in \Def(\MM)$ for all $b \in I$, then 
  its {\bf characteristic path} is the set 
 $\{s \in 2^{<I} : $ for all $a < \ell(s),  a \in A$ iff $ (s)_a = 0\}$.  The characteristic path of any such $A \subseteq I$ is a $2^{<I}$-path.

Suppose that $\MM$ is a model and $\Phi(x)$ is a set of 1-ary ${\mathcal L}(M)$-formulas. Let ${\mathcal C}(\Phi(x))$ be the set of 1-ary ${\mathcal L}(M)$-formulas $\varphi(x)$ such that $T_M \cup \Phi(x) \vdash \varphi(x)$. 
As usual, $\Phi(x)$ is {\bf consistent} if $x \neq x$ is not in ${\mathcal C}(\Phi(x))$, and 
$\Phi(x)$ is {\bf complete} if for every $1$-ary ${\mathcal L}(M)$-formula $\varphi(x)$, either 
$\varphi(x) \in {\mathcal C}(\Phi(x))$ or $\neg\varphi(x) \in {\mathcal C}(\Phi(x))$. 

If $\MM \preccurlyeq \NN \models \pa$, then the cut $\gcis(\MM,\NN)$ is defined in the Introduction.
  (It is possible to extend this definition to apply whenever there is $\KK$ such that $\MM, \NN \prec \KK$, but that will not be needed here.) If $I$ is a cut of $\MM$, $\MM \prec \NN$ and $b \in N$, then $b$ {\bf fills} $I$ if if $a < b < c$ whenever $a \in I$ and $c \in M \backslash I$. Thus,  $\NN$ fills $I$ iff some $b \in N$ does. Usually when saying that $\NN$ fills $I$, it will be that 
$I = \gcis(\MM,\NN)$. As defined in the introduction, if $\MM \prec_{\sf cf} \NN$, then $\NN$ is a filling extension of $\MM$ iff $\NN$ fills $\gcis(\MM, \NN)$, and $\NN$ is a non-filling extension otherwise.

A second-order model of arithmetic has the form $(\MM, {\mathfrak X})$, where ${\mathfrak X} \subseteq {\mathcal P}(M)$. Certain sets of second-order formulas, allowing parameters from $M$ and ${\mathfrak X}$, but with no  
second-order quantifiers, are the $\Sigma_n^0$ and $\Pi_n^0$, while each of $\Sigma_0^0, \Delta_0^0, 
\Pi_0^0$ is the set of all such formulas that are bounded. 
 For  $n < \omega$, let $\Sigma_n^0$-$\Def(\MM, {\mathfrak X})$  be the set of those $A \subseteq M$ that are definable in $(\MM, {\mathfrak X})$ by a 
$\Sigma_n^0$  formula; $\Pi_n^0$-$\Def(\MM, {\mathfrak X})$ is defined similarly; and 
$\Delta_n^0$-$\Def(\MM, {\mathfrak X}) = \Sigma_n^0$-$\Def(\MM, {\mathfrak X}) \cap \Pi_n^0$-$\Def(\MM, {\mathfrak X})$. A subset ${\mathfrak X}_0 \subseteq {\mathfrak X}$ {\bf generates} 
${\mathfrak X}$ if $\Delta_1^0$-$\Def(\MM, {\mathfrak X}_0) \supseteq {\mathfrak X}$. We say that 
${\mathfrak X}$ is $\kappa$-generated if it is generated by some ${\mathfrak X}_0$ such that 
$|{\mathfrak X}_0| \leq \kappa$. If ${\mathfrak X}$ is $\aleph_0$-generated,  then we say that it is {\bf countably generated}. 

Some  second-order theories will be important in this paper:

$$\rcas = \ISig_0^0 + {\sf exp} + \Delta^0_1\mbox{-}{\sf CA},$$
$$
\wkls = \rcas + {\sf WKL},
$$
$$\rca = \rcas + \ISig^0_1,$$
$$\aca = \rca +  \Sigma^0_1\mbox{-}{\sf CA}.$$

Fix a cut $I$ 
 of a model $\MM \models \pa$. Recall the following properties that $I$ may have: $I$ is {\bf semiregular} if for every $a,b \in M$ such that $ \ell(a) \geq b \in I$, there is $c \in I$ such that $(a)_i < c$ whenever $i < b$ and $(a)_i \in I$; $I$ is {\bf regular} if for every $a \in M$ such that $\ell(a) > I$, if  there is $b \in I$ such that $(a)_i < b$ for every $i \in I$, then  there is $c \in I$ such that $\{i \in I : (a)_i = c\}$ is an unbounded subset of $I$; $I$ is {\bf strong} if for every $a \in M$ such that $\ell(a) > I$, there is $b \in M$ such that for every $i \in I$, $(a)_i \in I$ iff $(a)_i < b$. Equivalently, $I$ is semiregular iff $(\MM | I, \cod(\MM / I)) \models \ISig_1^0$; 
 $I$ is regular iff  $(\MM | I, \cod(\MM / I)) \models \BSig_2^0$  iff  $(\MM | I, \cod(\MM / I)) \models \IDel_2^0$; 
$I$ is strong iff 
 $(\MM | I, \cod(\MM / I)) \models \aca$.

We end this section with  two propositions.

\bigskip

{\sc Proposition 1.1}: {\em Suppose that $\MM  \models \pa$ and $I$ is  an exponentially closed cut of $\MM$. Then$:$
\begin{itemize}

\item[(1)] $ (\MM | I , \cod(\MM / I)) \models \wkls$.

\item[(2)] If $\kappa$ is an infinite cardinal and $\MM$ is $\kappa$-generated over $I$, then 
$\cod(\MM /I)$ is $\kappa$-generated.

\end{itemize}}

\bigskip

{\it Proof.}  (1) This  is well known, being, for example,  a consequence of \cite[Th.\@ 4.8]{ss}. 
\smallskip

(2) Let $G \subseteq M$ be such that $G$ generates $\MM$ over $I$ and $|G| \leq \kappa$. 
We assume that $G$ is closed under the pairing function.
Let ${\mathfrak X} = \{A \cap I : A \in \Def(\MM)$ is definable by an ${\mathcal L}(G)$-formula$\}$.
Clearly, ${\mathfrak X} \subseteq \cod(\MM / I)$ and $|{\mathfrak X}| \leq \kappa$. We will show that 
${\mathfrak X}$ generates $\cod(\MM / I)$. 

Let $A \in \cod(\MM / I)$. Let $B \in \Def(\MM)$ be such that $A = B \cap I$. Let $\varphi(x,y)$ be an ${\mathcal L}(G)$-formula and $b \in I$ be such that $\varphi(x,b)$ defines $B$. Let $C \in \Def(\MM)$ be defined by $\varphi(x,y)$. Then $C \cap I \in \cod(\MM / I)$ and $A = (C)_b$.  \qed

\bigskip

{\sc Proposition 1.2}: {\em Suppose that   $\NN \succ_{\sf cf} \MM \models \pa$ and that $I = \gcis(\MM, \NN)$. Then$:$

\begin{itemize} 

\item[(1)]  $I$ is closed under multiplication. 

\end{itemize}
If, in addition, 
$\NN$ fills $I$, then$:$ 

\begin{itemize}

\item[(2)]  $I$ is regular. 

\item[(3)] $\Delta_2^0$-$\Def(\MM | I, \cod(\MM/I)) \subseteq \cod(\NN / I)$.
\end{itemize}}

\bigskip

{\it Proof}. The definition of $\gcis$ was introduced by Smory\'{n}ski \cite{smo}, where the easy proof \cite[Lemma~2.3]{smo}  of 
$(1)$ is given. It is also proved there that if $\NN$ fills $I$, then   $I$ is exponentially closed. The stronger conclusion of $(2)$ that  $(\MM | I, \cod(\MM / I)) \models \BSig^0_2$ is proved  in \cite[Th.\@ 3]{kp}. 
We prove (3).

Let $c \in N$ fill $I$.  Let $A \in \Delta_2^0$-$\Def(\MM | I, \cod(\MM/I))$. Then, there is $B \in \Def(\MM / I)$ such that
\begin{eqnarray}
A & = & \{a \in I : (\MM | I, \cod(\MM / I)) \models \forall x \exists y>x[\langle a,y \rangle \in B]\} \nonumber \\
& = & \{ a \in I :  (\MM | I, \cod(\MM / I)) \models \exists x \forall y>x[\langle a,y \rangle \in B]\}. \nonumber 
\end{eqnarray}
Let $b \in 2^{<M}$ be such that $\ell(b) >I$ and for every $i \in I$, $(b)_i = 0$ iff $i \in B$. One easily verifies that for any $a \in I$, $a \in A$ iff $\NN \models (b)_{\langle a,c \rangle} = 0$. Thus, $A \in \cod(\NN / I)$.
\qed

\bigskip

It should be noted that both $(2)$ and $(3)$  are implied by

\begin{itemize}

\item[(4)] $\Sigma_1^0$-$\Def(\MM | I, \cod(\MM/I)) \subseteq \cod(\NN / I)$

\end{itemize}
under the conditions of Proposition~1.2 that  $\NN \succ_{\sf cf} \MM \models \pa$ and  $I = \gcis(\MM, \NN)$. 
Clearly, $(4) \Longrightarrow (3)$. For (2), assume that $I$ is not regular and that $c \in M$ and $d \in I$ are such that $\ell(c) \geq d$ and $\{(c)_i : i \leq d\} \cap I$ is unbounded in $I$.  Let 
$X$ be the set of those $x <d$ such that $\{(c)_i : i \leq x\} \cap I$ is bounded. 
Then $X \in \Sigma_1^0$-$\Def(\MM | I, \cod(\MM/I)) \subseteq \cod(\NN / I)$. But $X \subsetneq I$ and $X$ is a cut of $\NN$, contradicting that $(\NN | I, \cod(\NN / I)) \models  \ISig_0^0$.

\smallskip

 The rest of this paper is devoted to obtaining various converses to Propositions~1.1 and~1.2 under the assumptions that $\MM$ is countably generated over $I$ and $I$ is an exponentially closed cut. 
  
 \bigskip


{\bf \S2.\@ The Core Theorem and its Corollaries.} This section contains one theorem and several  of its corollaries. For a model $\MM \models \pa$ that is countably generated over an exponentially closed cut $I$, the theorem (almost) characterizes those ${\mathfrak X}$ for which there is a finitely generated extension $\NN \succ_{\sf cf} \MM$ such that $I = \gcis(\MM, \NN)$ and ${\mathfrak X} = \cod(\NN / I)$. 

\bigskip

{\sc Theorem 2.1}:  {\em Suppose that $\MM \models \pa$,  $\MM$ is countably generated over the exponentially closed cut $I$ and $\cod(\MM / I) \subseteq {\mathfrak X} \subseteq {\mathcal P}(I)$.  
The following are equivalent$:$

\begin{itemize}

\item [(1)] $(\MM | I, {\mathfrak X}) \models \wkls$ and 
${\mathfrak X}$ is countably generated. 

\item[(2)] 
 There is a countably generated extension $\NN \succ_{\sf cf} \MM$ such that $\cod(\NN / I) = {\mathfrak X}$.

\item[(3)] 
 There is a  finitely generated extension $\NN \succ_{\sf cf} \MM$ such that $\cod(\NN / I) = {\mathfrak X}$. 
 
 \end{itemize}}

 \bigskip

 We introduce some terminology that will be used not only in the proof of this theorem but also in other proofs  in this paper.  
  \bigskip

{\sc Definition 2.2}:    Suppose that $\MM \models \pa$, $I$ is  a  cut of $\MM$,  $\cod(\MM / I)$ $ \subseteq {\mathfrak X} \subseteq {\mathcal P}(I)$ and $(\MM | I, {\mathfrak X}) \models \rcas$. We say that a set $\Phi(x)$ of $1$-ary ${\mathcal L}(M)$-formulas is {\bf allowable}   ({\bf for} $I$ {\bf and} ${\mathfrak X}$) if, for some $m,n < \omega$,  there are  $(1+n)$-ary ${\mathcal L}(M)$-formulas 
  $\theta_0(x,\overline u),$ $ \theta_1(x,\overline u), \ldots, \theta_{m-1}(x,\overline u)$ and  there are  $A_0,A_1, \ldots, A_{n-1} \in {\mathfrak X}$  such that     $$
  \Phi(x) =     \bigcup_{i < m}\{\theta_i(x,a_0,a_1, \ldots,a_{n-1}) : a_j \in A_j {\mbox{ for }} j < n\}.
  $$  
  
  Note that there is no requirement that an allowable set be consistent. When considering an allowable $\Phi(x)$, we are really more interested in ${\mathcal C}(\Phi(x))$ than  in $\Phi(x)$ itself.
  
  If we wish to emphasize the $m,n$ in Definition~2.2, then we say that $\Phi(x)$ is $(m,n)$-allowable. However, the following lemma shows that there is almost no need to do so.

  \bigskip

 {\sc Lemma 2.3}: {\em If $\Phi(x)$ is allowable, then there is a $(1,1)$-allowable $\Phi'(x)$ such that 
 ${\mathcal C}(\Phi'(x)) = {\mathcal C}(\Phi(x))$.}
 
  \bigskip
 
 {\it Proof}. Suppose $\MM \models \pa$ and that $\Phi(x)$ is $(m,n)$-allowable for $I$ and ${\mathfrak X}$ as in Definition~2.2. Let $A = A_0 \times A_1 \times \cdots \times A_{n-1} \in {\mathfrak X}$ and let 
 $\theta'(x,u)$ be the formula $\bigwedge_{i<m}\exists u_0,u_1, \ldots,u_{n-1}[u = \langle u_0,u_1, \ldots,u_{n-1} \rangle \wedge \theta_i( x, \overline u)]$. Let $\Phi'(x) = \{\theta'(x,a) : a \in A\}$. Then $\Phi'(x)$ is $(1,1)$-allowable and ${\mathcal C}(\Phi'(x)) = {\mathcal C}(\Phi(x))$. \qed

 \bigskip

It is easily verified that if $\Phi(x), \Phi'(x)$ are allowable sets and ${\mathcal C}(\Phi(x)) \subseteq  {\mathcal C}(\Phi'(x))$, then for any  $\Phi_0(x) \subseteq \Phi(x)$ that is definable (meaning definable in $\MM$), there is a definable $\Phi_0'(x) \subseteq \Phi'(x)$ such that if $X$ and $X'$ are defined by $\bigwedge \Phi_0(x)$ and $\bigwedge \Phi_0'(x)$, respectively, then $X' \subseteq X$.  (

\bigskip

The next definition  restricts the $(1,1)$-allowable sets even more. 

\bigskip

{\sc Definition} 2.4: (1) The $2$-ary ${\mathcal L}(M)$-formula 
  $\theta(x,u)$ is a {\bf tree-formula} if the three sentences
  $$
  \forall x,u[\theta(x,u) \into u \in 2^{<M}],
  $$
 $$
  \forall x,u,v[\theta(x,v) \wedge u \lhdeq v \into \theta(x,u)],
 $$
 $$
  \forall x,u,v[\big(\theta(x,u) \wedge \theta(x,v)\big) \into (u\lhdeq v \vee v \lhdeq u)]
 $$
 are in $T_M$. 
 
 (2) An allowable set $\Phi(x)$ is {\bf tree-based} if there is a tree-formula $\theta(x,u)$ and a $2^{<I}$-path $A \in {\mathfrak X}$ such that 
 $\Phi(x) = \{\theta(x,a) : a \in A\}$. 
 
 \bigskip
 
 Tree-based allowable sets are $(1,1)$-allowable.   It is not the case that for every allowable $\Phi(x)$ there is a tree-based allowable $\Phi'(x)$ such that ${\mathcal C}(\Phi'(x)) = {\mathcal C}(\Phi(x))$ (but see Lemma~2.1.3).
  
 \bigskip

We make one more definition.

\bigskip
 
{\sc Definition 2.5}: 
 If $\Phi(x)$ is an allowable  set for $I$ and ${\mathfrak X}$, then we say that 
$\varphi(x,v)$ {\bf represents} $A$ in $\Phi(x)$ if  $A \subseteq I$, $\varphi(x,v)$ is a $2$-ary ${\mathcal L}(M)$-formula and for all $a \in I$, 
$$
 \varphi(x,a) \in {\mathcal C}(\Phi(x)) \Longleftrightarrow   a \in A \Longleftrightarrow \neg\varphi(x,a) \not\in {\mathcal C}(\Phi(x)).
$$
If  some $\varphi(x,v)$  represents $A$ in $\Phi(x)$, then $A$ is {\bf represented} in $\Phi(x)$.

\bigskip

With these  definitions out of the way, we  return to the the proof of Theorem~2.1.

\bigskip

{\it Proof of  Theorem 2.1.} We fix $\MM$, $I$ and ${\mathfrak X}$ as in  Theorem~2.1.  The implication $(3) \Longrightarrow (2)$  is trivial, and   $(2) \Longrightarrow (1)$  is an immediate consequence of Proposition~1.1. The rest of this proof is devoted to  proving $(1) \Longrightarrow (3)$.

\smallskip

 Fix some $e \in M \backslash I$.  An allowable $\Phi(x)$   is $e$-{\bf big} if there is $d > I$ such that whenever $\Phi_0(x) \subseteq \Phi(x)$ is $\MM$-definable and $X = \{a <e : \MM \models \bigwedge \Phi_0(a)\}$,  then $\MM \models |X| \geq d$. Any such $d$ will be referred to as a {\bf bound} for $\Phi(x)$. Obviously, every 
 $e$-big allowable $\Phi(x)$ is consistent. The allowable set $\varnothing$ is $e$-big.

 \bigskip

 We will need the following two lemmas. 

\bigskip

{\sc Lemma 2.1.1}: {\em Suppose that  $\Phi(x)$ is an $e$-big allowable set. 
If $A \in {\mathfrak X}$, then there is   an $e$-big allowable $\Phi'(x) \supseteq \Phi(x)$  such that  $A$ is represented in $\Phi'(x)$.}

\bigskip 

{\sc Lemma 2.1.2}: {\em Suppose that  $\Phi(x)$ is an $e$-big allowable set. 
 If $\theta(x,u)$ is a $2$-ary ${\mathcal L}(M)$-formula, then there are  $A \in {\mathfrak X}$ and an $e$-big allowable $\Phi'(x) \supseteq \Phi(x)$  such that  $\theta(x,u)$ represents $A$ in $\Phi'(x)$.}

\bigskip
 
  Before proving these  lemmas, we will see how they are used to prove $(1) \Longrightarrow (3)$.  By the two countability conditions, we let $G \subseteq M \backslash I$ be a countable set that generates $\MM$ over $I$ and let ${\mathfrak X}_0 \subseteq {\mathfrak X}$ be a countable set that generates ${\mathfrak X}$. We can assume that $e \in G$ and also that $G$ is closed under  the pairing function. Let $\Phi_0(x) = \varnothing$, which is an $e$-big allowable set. Using that $G$ and ${\mathfrak X}_0$ are countable, we easily construct an  increasing sequence $\Phi_0(x) \subseteq \Phi_1(x) \subseteq \Phi_2(x) \subseteq \cdots$ 
  of $e$-big allowable sets  such that:

  \begin{itemize}

  \item[(S1)] For each $A \in {\mathfrak X}_0$,  there is $n < \omega$ such that   $A$ is represented in $\Phi_n(x)$.

\item[(S2)]  For each $2$-ary ${\mathcal L}(G)$-formula 
$\varphi(x,v)$, there are $A \in {\mathfrak X}$ and $n < \omega$ such that $\varphi(x,v)$ represents $A$ in $\Phi_n(x)$.

\end{itemize}

 Lemma~2.1.1 is used to get (S1), and Lemma~2.1.2 to get (S2).

Having this sequence, let $\Phi(x) =   \bigcup\{\Phi_n : n < \omega\}$. Clearly, $\Phi(x)$ is a consistent set of $1$-ary ${\mathcal L}(M)$-formulas since each $\Phi_n(x)$ is. We claim that $\Phi(x)$ is complete. To prove this claim, let $\psi(x)$ be  any $1$-ary ${\mathcal L}(M)$-formula. Then $\psi(x) = \varphi(x,a)$, 
where $\varphi(x,v)$ 
is an ${\mathcal L}(G)$-formula and $a \in I$. For this formula $\varphi(x,v)$, let $A$ and $n$ be as in (S2). If $a \in A$, then $\varphi(x,a) \in {\mathcal C}(\Phi_n(x))$, and if $a \not\in A$, then $\neg\varphi(x,a) \in {\mathcal C}(\Phi_n(x))$.
It then follows that either $\psi(x) \in {\mathcal C}(\Phi(x))$ or $\neg\psi(x) \in {\mathcal C}(\Phi(x))$, thereby proving the claim. A consequence of the completeness  is that  the formula $x < e$ is in ${\mathcal C}(\Phi(x))$.

Thus, ${\mathcal C}(\Phi(x))$ is a complete type over $\MM$. 
Let $\NN \succ \MM$ be an elementary extension generated by an element $c$ realizing this complete type. 
Obviously, $\NN$ is a finitely generated extension and $c < e$, so $\NN$ is a cofinal extension of $\MM$. We show that $\cod(\NN / I) =  {\mathfrak X}$. 

Suppose that $B \in {\mathfrak X}$. 
Since ${\mathfrak X}_0$ generates ${\mathfrak X}$, there are $A \in {\mathfrak X}_0$ and $a \in I$ such that $(A)_a = B$. For this $A$, let  $n < \omega$ be as in~(S1),    and then let $\varphi(x,v)$ represent $A$ in $\Phi_n(x)$. Let $X \in \Def(\NN)$ be defined by $\varphi(c,v)$. Then $A = X \cap I \in \cod(\NN / I)$. 
Then $B = \{b \in I : \NN \models \varphi(c,\langle a,b \rangle)\}$, so that $B \in {\cod(\NN / I)}$. 

Conversely, suppose that $B \in \cod(\NN / I)$. Let $Z \in \Def(\NN)$ be such that $B = Z \cap I$. 
Let the ${\mathcal L}(G)$-formula $\psi(x,y,z)$  and $a \in I$ be such that $\psi(c,a,z)$ defines $Z$ in $\NN$. Let $\varphi(x,v)$ be the ${\mathcal L}(G)$-formula $\exists y,z[\psi(x,y,z) \wedge v = \langle y,z \rangle]$. For this formula $\varphi(x,v)$, let $A \in {\mathfrak X}$ and $n < \omega$ be as in~(S2). Then $B = (A)_a$, so $B \in {\mathfrak X}$.

 Thus, $\NN$ is as  in (3) of Theorem~2.1. Incidentally, $I \subseteq \gcis(\MM,\NN)$ since $\cod(\NN / I) =  {\mathfrak X}$ and $(\MM | I, \mathfrak{X}) \models I\Sigma_0^0$.
 
 Modulo the proofs of Lemmas~2.1.1 and~2.1.2, we have proved Theorem~2.1.
We now turn to proving these two lemmas. (Neither of these lemmas need the countability conditions:   $\MM$ does not need to be countably generated over $I$ nor does  ${\mathfrak X}$ need  to be countably generated.) We begin with  a lemma that advantageously enables us to replace an allowable
set with a tree-based one, the advantage being that tree-based allowable sets are easier to work with than are arbitrary allowable sets.  For, if $\Phi(x) = \{\theta(x,a) : a \in A\}$ is tree-based and bounded by $d > I$, then $\Phi(x)$ is $e$-big iff there is $d > I$ such that for any $a \in A$, if $\theta(x,a) \wedge x<e$ defines $X_a$, then $\MM \models |X_a| \geq d$. Also, if $\Phi(x)$ is $e$-big (or even just consistent), then $\theta(x,u)$ represents $A$ in $\Phi(x)$.

 \bigskip
 
 {\sc Lemma 2.1.3}: {\em  If $\Phi(x)$ is $e$-big and allowable, then there is a tree-based, $e$-big allowable $\Phi'(x)$ such that 
 ${\mathcal C}(\Phi'(x)) \supseteq {\mathcal C}(\Phi(x))$.}
 
 \bigskip
 
 {\it Proof}. By Lemma~2.3 (and by the remark following the proof of that lemma), we assume that $\Phi(x)$ is $(1,1)$-allowable. Let $\Phi(x) = \{\theta(x,a) : a \in A\}$ be $(1,1)$-allowable and bounded by $d > I$. 
 Let $d_0$ be such that $I < 2^{d_0} < d$. Define 
 $\theta'(x,u)$ to be the formula 
 $$
 u \in 2^{<M} \wedge \forall i < \ell(u)[\theta(x,i) \leftrightarrow (u)_i = 0].
 $$ 
 It is easily seen that $\theta'(x,u)$ is a tree-formula.  Let $A' $ be the characteristic path of $A$, and let $\Phi'(x) = \{\theta'(x,s) : s \in A'\}$. 
 Clearly, $A' \in {\mathfrak X}$, so that $\Phi'(x)$ is an allowable set (and, hence,  tree-based).  One easily verifies that ${\mathcal C}(\Phi'(x)) \supseteq {\mathcal C}(\Phi(x))$ and that $\Phi'(x)$ is bounded by~$d_0$. \qed 

\bigskip

For the proofs of Lemmas~2.1.1 and~2.1.2, 
 we assume that $\Phi(x)$ is tree-based (by Lemma~2.1.3) and let 
$\Phi(x) = \{\theta_1(x,s) : s \in B\}$, where $\theta_1(x,u)$ is a  tree-formula and $B \in {\mathfrak X}$ is a $2^{<I}$-path.  For each $s \in M$, we let $X_s$ be the set defined by $\theta_1(x,s) \wedge x <e$. 
Let $d > I$ be a bound for $\Phi(x)$, and then let $d_0$ be such that $I < 2^{d_0} < d$. 

\smallskip
 
{\it Proof of Lemma 2.1.1}: Let $A  \in {\mathfrak X}$, and then let $A' \in {\mathfrak X}$ be the characteristic path of $A$.

 In $\MM$, we   try to define  by recursion two sequences $\langle b_i : i <e \rangle$ and 
$\langle y_{i,s} : i < e, \ s \in 2^{<M}, \ \ell(s) = \ell(b_i) \rangle$, where each $b_i \in 2^{<M}$ and each $y_{i,s} \in X_{b_i}$. Suppose that $i < e$ and that we already have $\langle b_j : j < i\rangle$ and $\langle y_{j,t} : j < i, \ t \in 2^{<M}, \ \ell(t) = \ell(b_j) \rangle$.
Let $b_i \in M$ be the least $b \in 2^{<M}$ such that 

\begin{itemize} 

\item $b \neq b_j$ for all $j < i$;

\item $\MM \models |\{X_b \backslash \{y_{j,t} : j< i, \ t \in 2^{<M}, \   b_j \lhd b, \ \ell(t) = \ell(b_j) \}| \geq 2^{\ell(b)}$.

\end{itemize}
 Then, in some uniformly definable way,  choose distinct $y_{i,s} \in X_{b_i} \backslash \{y_{j,t} : j< i, \ t \in 2^{<M}, \   b_j \lhd b, \ \ell(t) = \ell(b_j) \}$, where $s \in 2^{<M}$ and $\ell(s) = \ell(b_i)$.  If we can get $b_i$, then we will always be able to get the $y_{i,s}$'s. However,  it may be that for some $i< e$, there is no such $b_i$. Let $e_0 \leq e$ be the largest for which we can get $\langle b_i : i <e_0 \rangle$. Both of the sequences $\langle b_i : i <e_0 \rangle$ and 
$\langle y_{i,s} : i < e_0, \ s \in 2^{<M}, \ \ell(s) = \ell(b_i) \rangle$ are definable in $\MM$. Since $\Phi(x)$ is $e$-big, an easy calculation shows that for each $s \in B$, there is $i < e_0$ such that $b_i = s$. Thus, $e_0 > I$. 

Let $\theta_2(x,s,t)$ be the formula
$$
\exists i < e_0\exists t' \rhdeq t[b_i \rhdeq s \wedge \ell(t') = \ell(b_i) \wedge y_{i,t'} = x].
$$
Then  $\Phi'(x) = \Phi(x) \cup \{\theta_2(x,s,t) : s \in B, t \in A', \}$ is allowable and  
 $e$-big with a bound of $d_0$. 
 
 We show that the formula $ \theta_2(x,0,v)]$ 
  represents $A'$ in $\Phi'(x)$ according to Definition~2.5. 
  Consider any $a \in I$.  For one direction, suppose that $\neg\theta_2(x,0,a) \not\in {\mathcal C}(\Phi'(x))$.
  Thus, there are $s \in B$ and $t \in A'$ such that $\MM \models \forall x[\theta_2(x,s,t) \into \theta_2(x,0,a)]$. But then $0 \lhdeq s$ and $a \lhdeq t$, implying that $a \in A'$. 
  The other direction is trivial since if  $a \in A'$, then $\theta_2(x,0,a) \in \Phi'(x) \subseteq {\mathcal C}(\Phi'(x))$ since $0 \in B$.

  It easily follows that $A$ is represented in $\Phi'(x)$.\qed

\smallskip

{\it Proof of Lemma 2.1.2}:   Suppose that $\theta(x,u)$ is an ${\mathcal L}(M)$-formula. Let 
 $\theta_0(x,w,u)$ be the formula 
$$
 \theta_1(x,w) \wedge v \in 2^{<M} \wedge\ell(w) = \ell(u) \wedge \forall i < \ell(u)[\theta(x,i) \leftrightarrow (u)_i = 0].
 $$
 Let $T$ be the subtree of $2^{<I}$ consisting of those $t \in 2^{<I}$ such that for some $s \in B$,
 $\MM \models |\{x < e : \theta_0(x,s,t)| \geq d_0$ . Then $T \in {\mathfrak X}$ and is an unbounded subtree of $2^{<I}$, so there is a $T$-path $A' \in {\mathfrak X}$. Let 
 $$
 \Phi'(x) = \Phi(x) \cup\{\theta_0(x,s,t) : s \in B, t \in A'\}. 
 $$
 It is clear that $\Phi'(x)$ is an $e$-big allowable set having $d_0$ as a bound. 
 The formula $\exists w\theta_0(x,w,u)$ represents $A'$ in $\Phi'(x)$. Let $A \in {\mathfrak X}$ be the set whose characteristic path is $A'$. Then $\theta(x,u)$ represents $A$ in $\Phi'(x)$. 
  \qed
  
  \smallskip
 
 This completes the proof of Theorem~2.1. 
  \qqed 
  
  \bigskip
  
  In proving $(1) \Longrightarrow (3)$ of Theorem~2.1, if we did not try to satisfy (S2), then we would not need to use Lemma~2.1.2 and, thereby not need that $(\MM|I, {\mathfrak X}) \models {\sf WKL}$. Thus, we get the following corollary.
  
  \bigskip
  
  {\sc Corollary 2.6}:  {\em Suppose that $\MM \models \pa$,  $\MM$ is countably generated over the exponentially closed cut $I$ and $\cod(\MM / I) \subseteq {\mathfrak X} \subseteq {\mathcal P}(I)$.  
Suppose also that  $(\MM | I, {\mathfrak X}) \models \rcas$ and 
${\mathfrak X}$ is countably generated. Then 
there is a  finitely generated extension $\NN \succ_{\sf cf} \MM$ such that $\cod(\NN / I) \supseteq  {\mathfrak X}$.} \qed

\bigskip

We can get a corollary in the style of Harrington's Theorem that $\wkl$ is $\Pi_1^1$-conservative over ${\sf RCA}_0$.

\bigskip

{\sc Corollary 2.7}: {\em Suppose that $\MM \models \pa$,  $\MM$ is countably generated over the exponentially closed cut $I$, $\cod(\MM / I) \subseteq {\mathfrak X} \subseteq {\mathcal P}(I)$, ${\mathfrak X}$ is countably generated, and $(\MM | I, {\mathfrak X}) \models \rcas$. Then there is ${\mathfrak X}' $ such  that ${\mathfrak X} \subseteq {\mathfrak X}' \subseteq {\mathcal P}(I)$, ${\mathfrak X}'$ is countably generated and $(\MM | I, {\mathfrak X}') \models \wkls$.}

  \bigskip
  
  {\it Proof}. Get $\NN$ as in Corollary~2.6 and
let ${\mathfrak X}' = \cod(\NN /I)$.   Then $I \subseteq \gcis(\MM, \NN)$, so that $I$ is a cut of $\NN$. By Proposition~1.1, ${\mathfrak X}'$ is as required.   \qed

\bigskip
   
   We next present four  more corollaries of Theorem~2.1.
   
   \bigskip

{\sc Corollary 2.8}: (\cite[Th.\@ 4]{code}) {\em If $\MM\models \pa$ and 
$\Def(\MM) \subseteq {\mathfrak X} \subseteq {\mathcal P}(M)$, then the following are equivalent$:$

\begin{itemize}

\item[(1)]  $(\MM, {\mathfrak X}) \models \wkls$ and ${\mathfrak X}$ is  countably generated.  

\item[(2)] There is a countably generated extension $\NN \succ_{\sf end} \MM$ such that $\cod(\NN / M) = {\mathfrak X}$.

\item[(3)] There is a finitely generated extension $\NN \succ_{\sf end} \MM$ such that $\cod(\NN / M) = {\mathfrak X}$.

\end{itemize}}

\bigskip

{\it Proof}. Assume that $\MM$ and ${\mathfrak X}$ are as given. The implication $(3) \Longrightarrow (2)$  is trivial, and   $(2) \Longrightarrow (1)$  is an immediate consequence of Proposition~1.1. To prove $(1) \Longrightarrow (3)$, first use the MacDowell-Specker-Phillips Theorem for get a finitely generated, conservative  extension $\MM' \succ_{\sf end} \MM$, and then use $(1) \Longrightarrow (3)$ of Theorem~2.1  to get a finitely generated $\NN \succ_{\sf cf} \MM'$ such that $\cod(\NN / M) = {\mathfrak X}$. \qed

\bigskip

{\sc Corollary 2.9}: {\em Suppose that $\MM \models \pa$, $\MM$ 
is countably generated over the exponentially closed cut $I$, and  $ \cod(\MM / I) \subsetneq {\mathfrak X} \subseteq {\mathcal P}(I)$.   The following are equivalent$:$

\begin{itemize} 

\item[(1)] $(\MM|I, {\mathfrak X}) \models \wkls$ and ${\mathfrak X}$ is $\aleph_1$-generated.

\item [(2)] There is an ${\aleph}_1$-generated  $\NN \succ_{\sf cf} \MM$ such that $\gcis(\MM, \NN)= I$ and $\cod(\NN / I) = {\mathfrak X}$.

\end{itemize}}

\bigskip

{\it Proof}.  Let $\MM, I$ and ${\mathfrak X}$ be as given. The implication $(2) \Longrightarrow (1)$  is  an immediate consequence of Proposition~1.1. For the converse implication, suppose that $(1)$ holds. Let $\langle {\mathfrak X}_\alpha : \alpha < \omega_1 \rangle$ be a sequence 
of subsets of ${\mathfrak X}$ such that:

\begin{itemize}

\item ${\mathfrak X}_0 \supsetneq \cod(\MM / I)$;

\item if $\alpha < \beta < \omega_1$, then ${\mathfrak X}_\alpha \subseteq {\mathfrak X}_\beta$;

\item if $\alpha < \omega_1$, then ${\mathfrak X}_\alpha$ is countably generated and $(\MM | I, {\mathfrak X}_\alpha) \models  \wkls$;

\item ${\mathfrak X } = \bigcup\{{\mathfrak X}_\alpha : \alpha < \omega_1\}$.

\end{itemize}
By repeated applications of Theorem~2.1, get $\langle \MM_\alpha : \alpha < \omega_1 \rangle$ such that: 

\begin{itemize}

\item ${\MM}_0 = \MM$;

\item if $\alpha < \omega_1$, then $\MM_{\alpha+1}$ is a finitely generated, cofinal extension of  
$ \MM_\alpha$ such that $\cod(\MM_{\alpha+1} / I) = {\mathfrak X}_{\alpha+1}$;

\item if $\alpha < \omega_1$ is a limit ordinal, then $\MM_\alpha = \bigcup\{\MM_\beta : \beta < \alpha\}$.

\end{itemize}
Let $\NN = \bigcup \{\MM_\alpha : \alpha < \omega_1\}$. Since ${\mathfrak X}_0 \supsetneq \cod(\MM / I)$, it must be that $\gcis(\MM, \MM_1) = I$. We then have that $\gcis(\MM,\NN) = I$. Clearly, $\NN$ is as required. \qed

 \bigskip

{\sc Corollary 2.10}:  {\em If $\MM\models \pa$ and 
$\Def(\MM) \subseteq {\mathfrak X} \subseteq {\mathcal P}(M)$, then the following are equivalent$:$

\begin{itemize}

\item[(1)]  $(\MM, {\mathfrak X}) \models \wkls$ and ${\mathfrak X}$ is  $\aleph_1$-generated.  

\item[(2)] There is an $\aleph_1$-generated extension $\NN \succ_{\sf end} \MM$ such that \\ $\cod(\NN / M) = {\mathfrak X}$.

\end{itemize}}

\bigskip

{\it Proof}. The implication $(2) \Longrightarrow (1)$ follows from  Proposition~1.1. For the converse $(1) \Longrightarrow (2)$, if ${\mathfrak X} = \Def(\MM)$, then this is just the MacDowell-Specker-Phillips Theorem;  otherwise, this corollary follows from Corollary~2.9 in the same way that $(1) \Longrightarrow (3)$ of Corollary~2.8 follows from Theorem~2.1. \qed

\bigskip

Recall that for any $\MM \models \pa$, $\ssy(\MM)$ is a Scott set. Scott \cite{scott} proved a fundamental converse:  if ${\mathfrak X}$ is a  countable Scott set, then there is a prime model $\MM \models \pa$ such that 
$\ssy(\MM) = {\mathfrak X}$. 

\bigskip  

{\sc Corollary 2.11}: ({\it cf.} Knight \& Nadel \cite{KN}) {\em Suppose that  $\MM\models \pa$ is countable and nonstandard, ${\mathfrak X}$ is a Scott set,     ${\mathfrak X} \supseteq \ssy(\MM)$ and $|{\mathfrak X}| \leq \aleph_1$. 
  Then there is an extension $\NN \succ_{\sf cf} \MM$ such that $\ssy(\NN) = {\mathfrak X}$.}

\bigskip  

{\it Proof.} If ${\mathfrak X} = \ssy(\MM)$, apply $(1) \Longrightarrow (2)$ of Theorem~2.1. Otherwise, apply $(1) \Longrightarrow (2)$ of Corollary~2.10 with $I = \omega$. \qed

 \bigskip
 
 
 {\bf \S3.\@ The Exceptional Case.} Theorem~2.1 in the previous section   does not explicitly say anything  about whether or not $I = \gcis(\MM,\NN)$. 
  However, it must be that $I \subseteq \gcis(\MM, \NN)$. Clearly, if ${\mathfrak X} \neq \cod(\MM / I)$ in Theorem~2.1,  then  $I = \gcis(\MM, \NN)$.  But what happens otherwise? This question is answered in  Theorem~3.3. Furthermore, a consequence of  Lemma~3.1  is that the extension  in (3) of Theorem~3.3 is minimal.  The even stronger Corollary~3.4 follows from the implication $(1) \Longrightarrow (3)$ of Theorem~3.3.     Observe that in the hypothesis of the theorem, the cut $I$ is arbitrary so that it need not be closed under exponentiation.

  Since indiscernible types appear in Theorem~3.3, we make some brief remarks concerning them. 
  Recall some definitions that can be found, for example, in \cite[Chap.\@ 3]{ksbook}. Let $\MM \models \pa$ and let $p(x)$ be a $1$-type over $\MM$. For $1 \leq n < \omega$,  $p(x)$     is $n$-{\bf indiscernible} if it is  nonprincipal and whenever 
  $\theta(\overline x)$ is an $n$-ary ${\mathcal L}(M)$-formula, there is $\varphi(x) \in p(x)$ such that the sentence
  $$
  \forall \overline x, \overline y[\bigwedge_{i<n}\big(\varphi(x_i) \wedge \varphi(y_i)\big) \wedge \bigwedge_{i<n-1}\big(x_i < x_{i+1} \wedge y_i < y_{i+1}\big) \rightarrow \big(\theta(\overline x) \leftrightarrow \theta(\overline y)\big)].
  $$
is true in $\MM$.  
 If $1 \leq m < n < \omega$, then every $n$-indiscernible type is $m$-indiscernible. On the other hand, for every countable, nonstandard  $\MM\models \pa$, there are $n$-indiscernible types over $\MM$ that are not $(n+1)$-indiscernible. (See \cite[Th.\@ 3.1.10]{ksbook}.) 
  The type $p(x)$ is {\bf indiscernible} if it is $n$-indiscernible whenever $1\leq n< \omega$.

  For any set $X$, we let ${\mathbf B}_X$ be the Boolean lattice $({\mathcal P}(X), \cup, \cap)$. In particular, if $n < \omega$, then ${\mathbf B}_n$ is the finite Boolean lattice having exactly $n$ atoms. 
  
   The following lemma is probably well known.

  \bigskip
  
  {\sc Lemma 3.1:} {\em Suppose that $\MM \models \pa$ and that  $p(x)$ is an $(n+1)$-indiscernible type over $\MM$. If the extension $\NN \succ \MM$ is generated  by $a = \langle a_0,a_1, \ldots, a_{n-1} \rangle$, where the $a_i$'s are distinct elements of $N$ each realizing $p(x)$, then 
  $\Lt(\NN / \MM) \cong {\mathbf B}_n$.} 
  
   \bigskip
   
   {\it Proof sketch.} For $J \subseteq n$, let $\MM_J \subseteq \NN$ be generated by $\{a_j : j \in J\}$. 
   We claim that the function $J \mapsto \MM_J$ is an isomorphism from ${\mathbf B}_n$ onto $\Lt(\NN /\MM)$. To prove the claim,  it suffices to prove: if $f : M^n \into M$ is definable, then there are $\varphi(x) \in p(x)$ and $I \subseteq n$ such that whenever $\overline x, \overline y \in M^n$ are increasing $n$-tuples of elements satisfying $\varphi(x)$, then $f(\overline x) = f(\overline y)$ iff $x_i = y_i$ for all $i \in I$. 
   
Let $f : M^n \into M$ be definable in $\MM$.   For $i <n$, let $\theta_i(x_0,x_1, \ldots,x_n)$ be the 
formula 
$$
f(x_0,x_1, \ldots, x_{i-1},x_{i+1} ,x_n) \neq f(x_0,x_1, \ldots, x_i, x_{i+2}, \ldots, x_n).
$$
Let $\varphi(x) \in p(x)$ 
be such that it {\it forces} each $\theta_i(\overline x)$; that is, there is $I \subseteq n$ such that 
whenever $b_0 < b_1 < \cdots < b_n$ are elements of $M$ each satisfying $\varphi(x)$, then 
$\MM \models \theta_i(\overline b)$ iff $i \in I$. Then $I$ is as required in the lemma. \qed

\bigskip

   Letting $n = 1$ in the Lemma~3.1, we have that if the extension $\NN \succ \MM$ is generated  
   by an element realizing a $2$-indiscernible type $p(x)$ over $\MM$, then $\NN$ is a minimal extension of $\MM$ (or, equivalently, $p(x)$ is a {\bf selective} type over $\MM$).  The converse fails: for every countable nonstandard $\MM \models \pa$, there is a bounded selective type over $\MM$ that is not $2$-indiscernible over $\MM$.  (See \cite[Th.\@ 3.2.15]{ksbook}.) Incidentally, one proof of the MacDowell-Specker-Gaifman Theorem yields the following strengthening of its  improvement by Phillips: Every $\MM \models \pa$ has a conservative extension $\NN \succ_{\sf end} \MM$ generated by an element realizing an indiscernible type over $\MM$.   
   
   \bigskip
   
   Finite Ramsey's Theorem as formalized  in \pa\ is used to get indiscernible types. 
   When discussing Ramsey's Theorem, the following notation will be useful: If $X$ is a set linearly ordered by $<$ and $n < \omega$, we let $\langle X \rangle^n$ be the set of increasing $n$-tuples from $X$. There is a proof of Ramsey's Theorem with elementary bounds that can be formalized in \pa,  thereby yielding the following lemma.

\bigskip

{\sc Lemma 3.2}: {\em Suppose that $I$ is an exponentially closed cut of a model $\MM \models \pa$.
Suppose that $X \in \Def(\MM)$  is such that $\MM \models |X| = a$ for some $a > I$. Let $k < \omega$ and 
$\theta(\overline x,u)$ be a $(k+1)$-ary ${\mathcal L}(M)$-formula. Then there is $Y \in \Def(\MM)$ such that $Y \subseteq X$, $\MM \models |Y| = b$ for some $b > I$, and 
$$
\MM \models \theta(\overline c,e) \leftrightarrow \theta(\overline d,e)
$$
whenever $\overline c, \overline d \in \langle Y \rangle^k$ and $e \in I$.}
\qed

   \bigskip
 
 {\sc Theorem 3.3:} {\em Suppose that $\MM \models \pa$ and $\MM$ is countably generated over the cut $I$. The following are equivalent$:$
 
 \begin{itemize}
 
 \item[(1)] $\dcf(I) = \aleph_0$ and $I$ is exponentially closed. 
 
 \item[(2)] There is a countably generated extension $\NN \succ_{\sf cf} \MM$ such that  $\gcis(\MM, \NN) = I$ and $\cod(\NN / I) = \cod(\MM / I)$.
 
 \item[(3)]  There is a non-filling extension $\NN \succ_{\sf cf} \MM$ that is generated by an element realizing an indiscernible type over $\MM$ such that  $\gcis(\MM, \NN) = I$ and $\cod(\NN / I) = \cod(\MM / I)$.

 \end{itemize}}

    \bigskip
 
  {\it Proof}. Let $\MM$ and $I$ be  as in the theorem. Obviously,   $(3) \Longrightarrow (2)$, so it suffices to prove  $(2) \Longrightarrow (1)$ and $(1) \Longrightarrow (3)$. 
  
  \smallskip
  
  $(2) \Longrightarrow (1)$: Let $\NN$ be as in (2). 
   
   To prove that $I$ is exponentially closed, suppose, for a contradiction, that $a \in I < 2^a \in M$.
   Since $\gcis(\MM, \NN) = I$, we let $b \in N \backslash M$ be such that $I < b < 2^a$.  Let $B$ be the set of those $i  \in N$  for which the $i$-th digit in the binary expansion of $b$ is $1$.  Then $B \in \Def(\NN)$ and $B \subseteq I$, so $B \in  \Def(\MM)$. Let $c \in M$ be such that $B$ is the set of $i \in M$ for which the $i$-th digit in the binary expansion of $c$ is $1$ iff $i = B$.   But then $b = c$, so $b \in M$, which is a contradiction. 
   
   Next, we prove that $\dcf(I) = \aleph_0$. 
   For a contradiction, suppose that $\dcf(I) > \aleph_0$. (When we write $\dcf(I) > \aleph_0$, it is ambiguous as to whether we mean $\dcf^\MM(I) > \aleph_0$ or $\dcf^\NN(I) > \aleph_0$. But it makes no difference since $\dcf^\MM(I) = \dcf^\NN(I)$.) Since $\NN$ is countably generated over $M$ and $\MM$ is countably generated over $I$, then $\NN$ is countably generated over $I$. Let $G \subseteq N$ be a countable set that generates $\NN$ over $I$. Since $\dcf(I) > \aleph_0 = |G|$ and there are only countably many Skolem ${\mathcal L}(G)$-terms, there is some such term  $t(x)$ such that for arbitrarily small $b \in N \backslash M$, there is $i \in I$ 
   for which $\NN \models t(i) = b$. Since $\cod(\NN / I) = \cod(M/I)$, there is $a \in M$ such that $\ell(a) > I$ and  whenever $i,j \in I$, then $\MM \models c_{\langle i,j \rangle} = 0$ iff $\NN \models t(i) = j$. One easily shows that  formula 
   $$
   \forall u,v < x[(a)_{\langle u,v \rangle} = 0  \into t(u) = v]
   $$
   defines $I$ in $\NN$. But $I \not\in \Def(\NN)$. 
   
 \smallskip

  $(1) \Longrightarrow (3)$: We will first prove the weaker implication $(1) \Longrightarrow (3')$, where 
  $(3')$ is the weakening of (3) in which $\NN$ is required to be only a finitely generated extension  instead of being  generated by an element realizing an indiscernible type over $\MM$. 
  
Let $I$ be as in (1). Let $d_0 > d_1 > d_2 > \cdots$ be a decreasing sequence of elements of $M$  converging to $I$. Let $G \subseteq M$ be a  countable set that generates $\MM$ over $I$.  (We arrange that $\varnothing \neq G \subseteq M \backslash I$.) We will obtain a decreasing sequence $X_0 \supseteq X_1 \supseteq X_2 \supseteq \cdots$ of bounded sets in $\Def(\MM)$ such that:
 
 \begin{itemize}
 
  \item[(T1)] For each $n < \omega$ there is $i < \omega$ such that $\MM \models d_i \leq |X_n| \leq d_n$.
  
  \item[(T2)] For every  ${\mathcal L}(G)$-formula $\varphi(x,u)$ there is $n < \omega$ such that 
  for every $a \in I$,
  $$
  \hspace{20pt} \MM \models \forall x,y \in X_n[\varphi(x,a) \longleftrightarrow \varphi(y,a)].
 $$
 
   \item[(T3)] If $I$ is strong, then for every Skolem ${\mathcal L}(G)$-term $t(x,u)$ there are $n < \omega$ and $a,b \in M$ such that $ I < b \leq \ell(a)$, $(a)_u \in I$ for every $u \in I$,  and 
 $$
\hspace{20pt}  \MM \models \forall u < b[\big(\forall x \in X_n(t(x,u) \leq (a)_u)\big) \vee \big(\forall x \in X_n(t(x,u) \geq b)\big)].
 $$

  \end{itemize}
  
  Suppose that we have such a sequence. Let $\Phi(x)$ be the set of $1$-ary ${\mathcal L}(M)$-formulas 
  $\varphi(x)$ that define a superset of some $X_n$. Obviously, (T1) implies that $\Phi(x)$ is consistent, and  (T2) implies that $\Phi(x)$ is  complete. Let $\NN$ be generated over $\MM$ be an element $c$  realizing $\Phi(x)$. Clearly, $\NN$ is a finitely generated, cofinal extension of $\MM$. It also follows from (T2) that $\cod(\NN / I) \subseteq \cod(\MM / I)$, so that, in fact, $\cod(\NN / I) = \cod(\MM / I)$, implying that $I \subseteq \gcis(\MM,\NN)$.  Since the sequence $d_0 > d_1 > d_2 > \cdots$ converges to $I$, we get from (T1) that $I \supseteq \gcis(\MM, \NN)$, so that $I = \gcis(\MM,\NN)$.  
 
  We show that $\NN$ is  non-filling. There are two cases that depend on whether or not $I$ is strong.
  
  \smallskip
  
  {\em $I$ is not strong}: Let $e \in M$ be such that $\ell(e) > I$ and $I = \sup(\{(e)_i : i \in I\} \cap I) = \inf(\{(e)_i : i \in I\} \backslash I)$. For a contradiction, suppose that $d \in N$ fills $I$. 
  Let $A = \{i \in I : \NN \models (e)_i < d\}$, so that  $A \in \Cod(\NN / I) = \Cod(\MM / I)$. Let $a \in 2^{<M}$ code $A$; that is, if $i \in I$, then $i \in A \Longleftrightarrow (a)_i = 0$. Then the formula
  $$
  \forall y,z \leq x[(a)_y = 0 \wedge (a)_z = 1 \into y < z]
  $$
  defines $I$, which is a contradiction. 
  
  \smallskip
  
  {\em $I$ is strong}: Consider a typical $d \in N \backslash M$, and let $t(x,u)$ be a Skolem ${\mathcal L}(G)$-term and $i \in I$ be such that $\NN \models t(c,i) = d$.
  Letting $a,b$ be as in (T3), we easily see that $\NN \models d < (a)_i \vee d > b$, so $d$ does not fill $I$. 
  
  \smallskip
  
  Thus, $\NN$ is a finitely generated
extension such that $I = \gcis(\MM,\NN)$ and $\cod(\NN / I) = \cod(\MM / I)$.

  \smallskip
  
  We next turn to constructing the sequence of $X_n$'s by recursion. Let $X_0 = [0,d_0-1]^\MM$.  Suppose that we have $X_n$ and that $d_j \leq |X_n| \leq d_n$.  
   
   Suppose we are at a stage at which we are trying to fulfill (T2). Let $\varphi(x,u)$ be the first (in some given enumeration of the 2-ary  ${\mathcal L}(G)$-formulas) that has not yet been taken care of. 
   Choose a sufficiently large $i < \omega$. Working in $\MM$, let $F : X_n \into {\mathcal P}([0,d_i])$ be such that for $x \in X_n$, $F(x) = \{a \leq d_i : \varphi(x,a)\}$. Let $Y \subseteq X_n$ be the largest such that $F$ is constant on $Y$. Since, $i$ was chosen to be large enough, we can let $X_{n+1} = Y$. 
   
   Next,  suppose that $I$ is strong and we are at a stage at which we are trying to fulfill (T3). Let $t(x,u)$ be the first (in some given enumeration of the 2-ary  Skolem ${\mathcal L}(G)$-terms) that has not yet been taken care of. Working in $\MM$, we define by recursion $Y_s \subseteq X_n$ and $a_s \in M$ for each  $s \in 2^{<M}$. Let $Y_0 = X_n$. Suppose that $Y_s$ has been defined, where $\ell(s) = k$.  Let $s_0,s_1 \rhd s$ be such that $\ell(s_0) = \ell(s_1) = k+ 1$, $(s_0)_k=0$ and $(s_1)_k = 1$. Let $a_s$ be the least $a$ such that if $Z = \{x \in Y_s : t(x,k) \leq a\}$, then $2|Z| \geq |Y_s|$. Then let $Y_{s_0} = Z$ 
  and $Y_{s_1} = \{x \in Y_s : t(x,k) \geq a_s\}$. 
  Since $I$ is strong and closed under exponentiation, we let $b > I$ be such that whenever $s \in 2^{<I}$, then $a_s \in I$ iff $a_s < b$. Let $B$ be the $2^{<I}$-path such that if $s \in B$ and $\ell(s) = k$, then 
  $a_s <b$ iff there is $t \in B$ such that $\ell(t) = k+1$ and $(t)_k = 0$. 
  
  To get the generator of $\NN$ to realize an indiscernible type,    add the following  property to those that the sequence of $X_n$'s should have.

\begin{itemize}

\item[(T4)] For every  $k < \omega$ and $(k+1)$-ary ${\mathcal L}(G)$-formula $\varphi(\overline x,u)$ there is $n < \omega$ such that 
  for every $a \in I$,
  $$
  \hspace{20pt} \MM \models \forall \overline x,\overline y \in \langle X_n \rangle^k[\varphi(\overline x,a) \longleftrightarrow \varphi(\overline y,a)].
 $$
 
 \end{itemize}
 Satisfying (T4) is done in the same way that (T2) is satisfied, using Lemma~3.2. Notice that (T4) subsumes (T2) by taking $k = 1$. \qed

  \bigskip
   
{\sc Corollary 3.4}: {\em   Suppose that $\MM \models \pa$, $\MM$  is countably generated over the exponentially closed cut $I$ and  $\dcf(I) = \aleph_0$. Then, for every cardinal $\kappa > 0$, there is a non-filling extension $\NN \succ_{\sf cf} \MM$ such that   
   $\cod(\NN / I) = \cod(\MM / I)$, $\Lt(\NN / \MM) \cong {\mathbf B}_\kappa$, and 
   $\gcis(\MM, \NN_0) = I$ whenever $\MM \prec \NN_0 \preccurlyeq \NN$.} 
   
   \bigskip
   
   {\it Proof}. Let $p(x)$ be the  type of a generator as in $(3)$ of Theorem~3.3. Now let $\NN \succ \MM$ be an extension generated by a set $X$ such that $|X| = \kappa$ and each element of $X$ realizes $p(x)$ over $\MM$. Using Lemma~3.1, one easily verifies that $\NN$ is as required. \qed

  \bigskip

Paris \& Mills \cite{pm} make the following definitions for a model $\MM \models \pa$ and an infinite cardinal $\kappa$:
$$ I^\MM_\kappa = \{x \in M : |[0,x]^\MM| \leq \kappa\},$$
$$ J^\MM(\kappa) = \inf\{|[0,a]^\MM| : a \in M \backslash I^\MM_\kappa\}.$$ 
If $I^\MM_\kappa \neq M$, then $I^\MM_\kappa$ is a cut closed under multiplication and $J^\MM(\kappa)$ is well defined. It is noted in \cite{pm} that if $J^\MM(\kappa) > 2^\kappa$, then 
$I^\MM_\kappa$ is closed under exponentiation.  

\bigskip

{\sc Corollary 3.5:} {\em Suppose $\MM \models \pa$,   
$\MM$ is countably generated over the exponentially closed cut $I$ and $\dcf(I) = \aleph_0$. 
Whenever $\kappa \geq \lambda = |I|$, there is a non-filling $\NN \succ_{\sf cf} \MM$ such that 
$I = \gcis(\MM, \NN) = I^\NN_\lambda$, $J^\NN(\lambda)= \kappa$ and $\cod(\MM/I) = \cod(\NN/I)$.} 

\bigskip

{\it Proof}. Let $\NN$ be as in Corollary~3.4. \qed

\bigskip

The previous corollary reduces to \cite[Th.\@~6]{pm} when $\MM$ is countable, although that theorem  does not explicitly  state that $\NN$ is a non-filling extension nor that $\cod(\NN/I) = \cod(\MM/I)$. 
This corollary is improved in Corollary~5.8.

\bigskip


 {\bf \S4.\@ Filling Extensions.} Recall from the introduction that if $\MM \prec_{\sf cf} \NN$, then $\NN$ is a filling extension of $\MM$ iff there is $b \in N$ such that $a < b < c$ whenever $a \in \gcis(\MM,\NN) < c \in M$. The next theorem characterizes when the extension $\NN \succ_{\sf cf} \MM$ in Theorem~2.1 can be filling. It will be handy to have the following definition: $\NN$ is an ${\mathfrak X}$-{\bf extension} of $\MM$ if $\NN \succ \MM \models \pa$ and  $\MM$ has a  cut $I$  such that  $I = \gcis(\MM, \NN)$ and ${\mathfrak X} = \cod(\NN / I)$.

\bigskip

 {\sc Theorem 4.1:} {\em Suppose that $\MM \models \pa$,  $\MM$ is countably generated over the cut $I$, $\cod(\MM / I) \subseteq {\mathfrak X} \subseteq {\mathcal P}(I)$ and $\MM$ has a countably generated ${\mathfrak X}$-extension.  The following are equivalent$:$  
 
 \begin{itemize}
 
 \item[(1)]    $\Delta_2^0$-$\Def(\MM | I, \cod(\MM/I)) \subseteq {\mathfrak X}$.

 \item[(2)]  There is a countably generated, filling extension $\NN \succ_{\sf cf} \MM$ such that  $\cod(\NN / I) = {\mathfrak X}$.
 
  \item[(3)]  There is an  extension $\NN \succ_{\sf cf} \MM$ such that  $\cod(\NN / I) = {\mathfrak X}$ 
  and $\NN$ is generated over $\MM$ by an element filling $I$.

 \end{itemize}}   
 
 \bigskip
 
 {\sc Remark}: Statement (3) is, on its face, stronger than just asserting that there is a finitely generated, filling extension $\NN \succ_{\sf cf} \MM$ such that  $\cod(\NN / I) = {\mathfrak X}$.
 
\bigskip

{\it Proof}. Let $\MM$, $I$ and ${\mathfrak X}$ be as in the theorem. The implication $(2) \Longrightarrow (1)$ follows from Proposition~1.2, and $(3) \Longrightarrow (2)$ is trivial. It remains to prove  $(1) \Longrightarrow (3)$. 

For the proof of $(1) \Longrightarrow (3)$, we say that an allowable set $\Phi(x)$ is $I$-{\bf big} if   whenever $\Phi_0(x) \subseteq \Phi(x)$ is $\MM$-definable and $X \in \Def(\MM)$ is defined by  $ \bigwedge \Phi_0(x)$,  then  $X \cap I$ is an unbounded subset of $I$.  We will need the following two lemmas that are analogous to Lemmas~2.1.1 and~2.1.2.

\bigskip

{\sc Lemma 4.1.1}: {\em Suppose that  $\Phi(x)$ is an $I$-big allowable set. 
If $A \in {\mathfrak X}$, then there is   an $I$-big allowable $\Phi'(x) \supseteq \Phi(x)$  such that  $A$ is representable in $\Phi'(x)$.}

\bigskip 

{\sc Lemma 4.1.2}: {\em Suppose that  $\Phi(x)$ is an $I$-big allowable set. 
 If $\theta(x,u)$ is a $2$-ary ${\mathcal L}(M)$-formula, then there are  $A \in {\mathfrak X}$ and an $I$-big allowable $\Phi'(x) \supseteq \Phi(x)$  such that  $\theta(x,u)$ represents $A$ in $\Phi'(x)$.}

\bigskip

Let $G$ and ${\mathfrak X}_0$ be as in the proof of Theorem~2.1 (except that we do not have that $e \in G$). 
Having these two lemmas, we construct an increasing sequence $\Phi_0(x) \subseteq \Phi_1(x) \subseteq \Phi_2(x) \subseteq \cdots$ of $I$-big allowable sets such that (S1) and (S2) from the proof of Theorem~2.1 hold. Let $\Phi(x) = \bigcup\{\Phi_n(x) : n < \omega\}$. Just as in the proof of Theorem~2.1, $\Phi(x)$ is a complete type over $\MM$, and the extension $\NN$ generated over $\MM$ by an element $c$ realizing $\Phi(x)$ is such that $\cod(\NN / I) = {\mathfrak X}$. Clearly, whenever $a \in I < b \in M$, then the formula $a < x < b$ is in 
${\mathcal C}(\Phi(x))$. Consequently,  $c$ fills $I$ and $I \subseteq \gcis(\MM,\NN)$, so that, in fact $I = \gcis(\MM,\NN)$ and  $\NN$ is a filling extension of $\MM$. 

Modulo the proofs of Lemmas~4.1.1 and~4.1.2, the proof of Theorem~4.1 is complete.  We now turn to proving these two lemmas. The following lemma is analogous to Lemma~2.1.3.

\bigskip
 
 {\sc Lemma 4.1.3}: {\em  If $\Phi(x)$ is $I$-big and allowable, then there is a tree-based, $I$-big allowable $\Phi'(x)$ such that 
 ${\mathcal C}(\Phi'(x)) \supseteq {\mathcal C}(\Phi(x))$.}
 
 \bigskip
 
 {\it Proof}.  By Lemma~2.3 (and by the remark following the proof of that lemma), we assume that $\Phi(x)$ is $(1,1)$-allowable.   Let $\Phi(x) = \{\theta(x,a) : a \in A\}$, where $A \in {\mathfrak X}$.
  Define 
 $\theta'(x,u)$ to be the formula 
 $$
 u \in 2^{<M} \wedge \forall i < \ell(u)[\theta(x,i) \leftrightarrow (u)_i = 0].
 $$ 
 It is easily seen that $\theta'(x,u)$ is a tree-formula.  Let $A' \in {\mathfrak X}$ be the characteristic path of $A$, and let $\Phi'(x) = \{\theta'(x,s) : s \in A'\}$. 
 Clearly,  $\Phi'(x)$ is an allowable set (and, hence,  tree-based).  One easily verifies that ${\mathcal C}(\Phi'(x)) \supseteq {\mathcal C}(\Phi(x))$ and that $\Phi'(x)$ is $I$-big. \qed 

\bigskip

For  the proofs of both Lemmas~4.1.1 and~4.1.2, we assume that $\Phi(x)$ is tree-based. Let 
$\Phi(x) = \{\theta_1(x,s) : s \in B\}$, where $\theta_1(x,u)$ is a tree-formula and $B \in {\mathfrak X}$ is a $2^{<I}$-path. For each $s \in M$, we let $X_s = \{x \in M : \MM \models \theta_1(x,s)\}$.

\smallskip

{\it Proof of Lemma~4.1.1}. This proof is  much like the proof of Lemma 2.1.1. Let $A \in {\mathfrak X}$,  and  then let $A' \in {\mathfrak X}$ be the characteristic path of $A$. 

In $\MM$, 
  we  try to inductively define two sequences $\langle b_i : i \in M \rangle$ and 
$\langle y_{i,s} : i \in M, \ s \in 2^{<M}, \ \ell(s) = \ell(b_i) \rangle$ as follows.
Suppose that $i \in M$ and that we already have $\langle b_j : j < i\rangle$ and $\langle y_{j,t} : j<i, \ t \in 2^{<M}, \ \ell(t) = \ell(b_j) \rangle$. Let $b_i$ be the least $b \in 2^{<M}$ such that:

\begin{itemize} 

\item $b > b_j$ for all $j < i$ such that $b_j \lhd b$.

\item $\MM \models |\{X_b \backslash \{y_{j,t} : j< i, \ t \in 2^{<M}, \   b_j \lhd b, \ \ell(t) = \ell(b_j) \}| \geq 2^{\ell(b)}$.

\end{itemize}

 Then, in some uniformly definable way,  choose $y_{i,s} \in X_{b_i} \backslash \{y_{j,t} : j<i, t \in 2^{<M},  \ b_j \lhd b_i, \ \ell(t) = \ell(b_j)\}$, where $s \in 2^{<M}$ and $\ell(s) = \ell(b_i)$, so that $s \mapsto y_{i,s}$ is one-to-one.   It may be that for some $i \in M$, there is no such $b_i$. If so, let $e_0 \in M$ be the largest for which we can get $\langle b_i : i <e_0 \rangle$ and $\langle y_{i,s} : i < e_0, \  s \in 2^{<M}, \ \ell(s) = \ell(b_i)\rangle$. (If not, then let $e_0 > I$ be arbitrary). Both of the sequences  $\langle b_i : i <e_0 \rangle$ and 
$\langle y_{i,s} :  i  < e_0, \ s \in 2^{<M}, \  \ell(s) = \ell(b_i)\rangle$ are definable in $\MM$. For each $s \in B$, there is $i < e_0$ such that $b_i = s$. Thus, $e_0 > I$. 

We continue just as in the proof of Lemma~2.1.1. 
Let $\theta(x,u)$ be the formula 
$$
 \exists i < e_0 \exists s[  u \lhdeq s \in 2^{<M} \wedge \ell(s) = \ell(b_i) \wedge x = y_{i,s}],
$$
and then let $\Phi'(x) = \Phi(x) \cup \{\theta(x,s) : s \in A'\}$. One easily sees that $\Phi'(x)$ is allowable and $I$-big. Just as in the proof of Lemma~2.1.1, $A$ is representable in $\Phi'(x)$. 
 \qed

\smallskip

{\it Proof of Lemma~4.1.2}. Suppose that $\theta(x,u)$ is  a $2$-ary ${\mathcal L}(M)$-formula. 
We will obtain an $I$-big allowable set $\Phi'(x) \supseteq \Phi(x)$ and $A \in {\mathfrak X}$ such that 
$\theta(x,u)$ represents $A$ in $\Phi'(x)$.

Let $D \in \Def(\MM)$ be defined by the formula
$$
v \in 2^{<M} \wedge \ell(u) = \ell(v) \wedge \theta_1(x,u) \wedge \forall w < \ell(v)[\theta(x,w) \leftrightarrow (v)_w = 0].
$$
Then $D \cap I \in \cod(\MM / I)$. Let $P = \{\langle u,v \rangle \in I : (\MM |I, \cod(\MM / I)) \models \exists x > \ell(u)[\langle x,u,v \rangle \in D\}$. Thus, $P \in \Sigma_1^0$-$\Def(\MM | I, \cod(\MM / I)) \subseteq {\mathfrak X}$. Let $T = \{v \in I :  (\MM |I, \cod(\MM / I)) \models \exists s \in B[\langle s,v \rangle \in P\}$.
Clearly, $T \in {\mathfrak X}$ and $T$ is an unbounded $I$-tree. Let $A' \in {\mathfrak X}$ be a $T$-path, 
and let $A \in {\mathfrak X}$ be such that $A'$ is the characteristic path of $A$. 
Then, $\Phi'(x) = \Phi(x) \cup \{\theta(x,a) : a \in A\} \cup \{\neg\theta(x,a) : a \in I \backslash A\}$ is an $I$-big allowable set in which $\theta(x,u)$ represents $A$. \qed

\smallskip

This completes the proof of Theorem~4.1. \qqed

\bigskip  

The notion of an extendible cut is implicit in \cite{kp}. If $I$ is a cut of a model $\MM$, then $I$ is {\bf extendible} if there is a filling extension $\NN \succ_{\sf cf} \MM$ such that $I = \GCIS(\MM,\NN)$. This definition was generalized in \cite{p}.

\bigskip

{\sc Definition 4.2}: Suppose that $\MM \models \pa$. We define, by recursion on $n < \omega$, when a cut $I$ is $n$-extendible in $\MM$.
Every cut $I$ of $\MM$ is $0$-{\bf extendible} in $\MM$. A cut $I$ is $(n+1)$-{\bf extendible} in $\MM$ if there is  a filling extension $\NN \succ_{\sf cf} \MM$ such that $I = \GCIS(\MM,\NN)$ and $I$ is an $n$-extendible cut of  $\NN$. 

\bigskip

The previous definition suggests the next one.

\bigskip

{\sc Definition 4.3}: The extension $\NN \succ_{\sf cf} \MM$ is $n$-{\bf filling} if there are $\MM = \NN_0 \prec_{\sf cf} \NN_1 \prec_{\sf cf} \dots \prec_{\sf cf} \NN_n = \NN$ such that for 
$i < n$, $\gcis(\NN_i, \NN_{i+1})$ $ = \gcis(\MM, \NN)$ and $\NN_{i+1}$ is a filling extension of $\NN_i$.

\bigskip

These two definitions are closely related by the next proposition, which is an immediate consequence of the definitions.

\bigskip

{\sc Proposition 4.4}: {\em Suppose that $I$ is a cut of $\MM \models \pa$ and $n < \omega$. Then, 
$I$ is $n$-extendible iff there is an $n$-filling extension $\NN \succ_{\sf cf} \MM$ such that $\gcis(\MM, \NN) = I$.} \qed

\bigskip

{\sc Corollary 4.5}: {\em Suppose that $1 \leq n < \omega$, $\MM \models \pa$, $\MM$ is countably generated over the exponentially closed cut $I$, $\cod(\MM / I) \subseteq {\mathfrak X} \subseteq {\mathcal P}(I)$ and $\MM$ has a countably generated ${\mathfrak X}$-extension. The following are equivalent$:$

\begin{itemize}

\item[(1)] $\Delta_{n+1}^0$-$\Def(\MM |I, \cod(\MM / I)) \subseteq {\mathfrak X}$.

\item[(2)] There is a countably generated, $n$-filling extension $\NN \succ_{\sf cf} \MM$ such that $\gcis(\MM, \NN) = I$ and $\cod(\MM / I) = {\mathfrak X}$.

\item[(3)] There is an $n$-filling extension $\NN \succ_{\sf cf} \MM$ such that $\gcis(\MM, \NN)$ $ = I$, 
 $\cod(\MM / I) = {\mathfrak X}$ and $\NN$ is generated over $\MM$ by an element filling $I$.

\end{itemize}}

\bigskip

{\it Proof}. The proof is by induction on $n$. The $n = 1$ case is exactly Theorem~4.1.  As an inductive hypothesis, suppose that $2 \leq m < \omega$ and that the corollary holds when $n = m-1$. Suppose that 
$n= m$ and that $\MM, I$ and ${\mathfrak X}$ are as given. The implication $(3) \Longrightarrow (2)$ is trivial. We prove $(2) \Longrightarrow (1)$ and $(1) \Longrightarrow (3)$.

\smallskip

$(2) \Longrightarrow (1)$: Let $\NN$ be as in (2). Then there is an $(n-1)$-filling extension $\NN_{n-1} \succ_{\sf cf} \MM$ such that $I = \gcis(\MM,\NN_{n-1})$ and $\NN \succ_{\sf cf} \NN_{n-1}$ is a filling extension. By the inductive hypothesis,  $\Delta_n^0$-$\Def(\MM | I, \cod(\MM / I)) \subseteq \cod(\NN_{n-1} / I)$ and, by Theorem~4.1, $\Delta_2^0$-$\Def(\MM | I, \cod(\NN_{n-1} / I)) \subseteq 
\cod(\NN /I)$. Therefore, $\Delta_{n+1}^0$-$\Def(\MM | I, \cod(\MM / I)) \subseteq \cod(\NN / I)$.

\smallskip

$(1) \Longrightarrow (3)$: Suppose that $(1)$ holds. There is ${\mathfrak X}_{n-1} \subseteq {\mathfrak X}$ such that $\Delta_{n}^0$-$\Def(\MM | I, \cod(\MM / I)) \subseteq {\mathfrak X}_{n-1}$ and 
$\Delta_2^0$-$\Def(\MM | I, {\mathfrak X}_{n-1}) \subseteq {\mathfrak X}$ and 
such that $(\MM | I, {\mathfrak X}_{n-1}) \models \wkls$ and ${\mathfrak X}_{n-1}$ is countably generated. 
By the inductive hypothesis, there is an $(n-1)$-filling ${\mathfrak X}_{n-1}$-extension $\NN_{n-1} \succ_{\sf cf} \MM$. By Theorem~5.6, let $\NN \succ_{\sf cf} \NN_{n-1}$ be a filling ${\mathfrak X}$-extension. Then $\NN$ is as required by (3). \qed

\bigskip

The following corollary  was  proved by Clote \cite{c} for countable $\MM$ after some earlier, 
partial progress by Kirby \& Paris \cite{kp} and Paris \cite{p}. 

\bigskip

{\sc Corollary 4.6}: {\em Suppose that $1 \leq n < \omega$, $\MM \models \pa$ and $\MM$ is countably generated over the  cut $I$.  The following are equivalent$:$

$(1)$ $I$ is $n$-extendible.

$(2)$  $(\MM| I,  \cod(\MM / I)) \models \BSig_{n+1}^0$.}

\bigskip

{\it Proof}.  Let $n, \MM$ and $I$ be as given.

\smallskip

$(1) \Longrightarrow (2)$: Suppose that $I$ is $n$-extendible. By Proposition~4.4, let $\NN \succ_{\sf cf} \MM$ be a finitely generated, $n$-filling extension such that $\gcis(\MM, \NN) = I$. Let ${\mathfrak X} = \cod(\NN / I)$. Then, $(\MM|I, {\mathfrak X}) \models \rcas$. By Corollary~4.5, $\Delta_{n+1}^0$-$\Def(\MM | I, \cod(\MM /I)) \subseteq {\mathfrak X}$, so  $(\MM| I,  \cod(\MM / I))$ $ \models \IDel_{n+1}^0$.
Therefore, $(\MM| I,  \cod(\MM / I)) \models \BSig_{n+1}^0$.

\smallskip

$(2) \Longrightarrow (1)$: Let ${\mathfrak X} = \Delta_{n+1}^0$-$\Def(\MM | I, \cod(\MM / I))$. Since 
$(\MM| I,$ $  \cod(\MM / I)) \models \BSig_{n+1}^0$ iff $(\MM| I,  \cod(\MM / I)) \models \IDel_{n+1}^0$, 
then $(\MM | I, {\mathfrak X}) \models \rca$. Apply Corollary~2.7 to get a countably generated ${\mathfrak X}' \supseteq {\mathfrak X}$. Corollary~4.5 implies that there is a finitely generated, $n$-filling extension $\NN \succ_{\sf cf} \MM$ such that $\gcis(\MM, \NN) = I$ and $\cod(\NN / I) = {\mathfrak X}'$. Thus, $I$ is $n$-extendable. \qed

\bigskip

The following definition is from Paris \cite{p890}.

\bigskip

{\sc Definition 4.7}: Suppose that $\MM \models \pa$. We define, by recursion on $n < \omega$, when a cut $I$ is $(n + \half)$-extendible in $\MM$.
A cut $I$ of $\MM$ is $\half$-{\bf extendible} in $\MM$ if $I$ is semiregular, and it is 
 $(n+1+ \half)$-{\bf extendible} in $\MM$ if there is  a filling extension $\NN \succ_{\sf cf} \MM$ 
such that $I = \GCIS(\MM,\NN)$ and $I$ is $(n+ \half)$-extendible in~$\NN$. 

\bigskip

The next definition and proposition are related to the previous definition in the same way that 
Definition~4.3 and Proposition~4.4 are related to Definition~4.2. 

\bigskip

{\sc Definition 4.8}: The extension $\NN \succ_{\sf cf} \MM$ is $(n+\half)$-{\bf filling} if $I$ is semiregular in $\MM$ and there are $\MM = \NN_0 \prec_{\sf cf} \NN_1 \prec_{\sf cf} \dots \prec_{\sf cf} \NN_n = \NN$ such that for 
$i < n$, $\gcis(\NN_i, \NN_{i+1}) = \gcis(\MM, \NN)$ and $\NN_{i+1}$ is a filling extension of $\NN_i$.

\bigskip

These two definitions are closely related by the next proposition, which is an immediate consequence of the definitions..

\bigskip

{\sc Proposition 4.9}: {\em Suppose that $I$ is a cut of $\MM \models \pa$ and $n < \omega$. Then, 
$I$ is $(n+ \half)$-extendible iff there is an $(n+ \half)$-filling extension $\NN \succ_{\sf cf} \MM$ such that $\gcis(\MM, \NN) = I$.} \qed

\bigskip

{\sc Corollary 4.10}: {\em Suppose that $1 \leq  n < \omega$, $\MM \models \pa$, $\MM$ is countably generated over the exponentially closed cut $I$, $\cod(\MM / I) \subseteq {\mathfrak X} \subseteq {\mathcal P}(I)$ and $\MM$ has a countably generated ${\mathfrak X}$-extension. The following are equivalent$:$

\begin{itemize}

\item[(1)] $\Delta_{n+1}^0$-$\Def(\MM |I, \cod(\MM / I)) \subseteq {\mathfrak X}$ and $I$ is semiregular. 

\item[(2)] There is a countably generated, $(n+ \half)$-filling extension $\NN \succ_{\sf cf} \MM$ such that $\gcis(\MM, \NN) = I$.

\item[(3)] There is an $(n+ \half)$-filling extension $\NN \succ_{\sf cf} \MM$ such that $\gcis(\MM, \NN)$ $ = I$ and $\NN$ is generated over $\MM$ by an element filling $I$.

\end{itemize}}

\bigskip

{\it Proof}. This follows from Corollary~4.5. \qed

\bigskip

The following corollary is analogous to Corollary~4.6. In a concluding remark in \cite{c}  its author stated that he had proved the following result for countable $\MM$. However, it seems that no published account has appeared.

\bigskip

{\sc Corollary 4.11}: {\em Suppose that $ n < \omega$, $\MM \models \pa$ and $\MM$ is countably generated over the cut $I$.  The following are equivalent$:$

$(1)$ $I$ is $(n+ \half)$-extendible.

$(2)$ 
 $(\MM | I, \cod(\MM / I)) \models \ISig_{n+1}^0$.} 

\bigskip

{\it Proof}. This proof is just like the proof of Corollary~4.7. It is helpful to observe that 
$\Pi_{n+1}^0$-$\Def(\MM | I, \cod(\MM / I)) = \Delta_{n+1}^0$-$\Def(\MM |I, \Sigma_1^0$-$\Def(\MM|I, \cod(\MM/I)))$.  \qed

\bigskip


 {\bf \S5.\@ Non-filling Extensions.} The primary purpose of this section is to characterize in Theorem~5.6 when the extension $\NN \succ_{\sf cf} \MM$ in Theorem 2.1 can be non-filling (except for the case that ${\mathfrak X} = \cod(\MM / I)$, which is taken care of by Theorem~3.3).  We begin  with  another definition. 

\bigskip

   {\sc Definition 5.1}: If $\kappa > 0$ is a cardinal,  then we say that 
   $(\MM, {\mathfrak X})$ is $\kappa$-{\bf strong} if every set in $\Pi_1^0$-$\Def(\MM, {\mathfrak X)}$ is the union of a at most $\kappa$ sets in $\Sigma_1^0$-$\Def(\MM, {\mathfrak X)}$.

   \bigskip

   Notice that if $0 < \kappa < \lambda$ and $(\MM, {\mathfrak X})$ is $\kappa$-strong, then  $(\MM, {\mathfrak X})$ is $\lambda$-strong. If $(\MM, {\mathfrak X}) \models \rca$, then  
   $(\MM, {\mathfrak X})$ is $1$-strong iff $(\MM, {\mathfrak X}) \models \aca$.
   This choice of terminology was suggested by the fact that if $\MM \models \pa$ and $I$ is an exponentially closed cut of $\MM$, then $(\MM|I, \cod(\MM / I))$ is $1$-strong iff $I$ is a strong cut. 
   
   Although not formally made in \cite{min}, this definition played a role there where the following is proved.
   
   \bigskip
   
   {\sc Theorem 5.2}: ({\it cf.\@} \cite[Th.\@ 3]{min})
   {\em Suppose that  $\MM \models \pa$, ${\mathfrak X} \subseteq {\mathcal P}(I)$ and there is a countably generated ${\mathfrak X}$-extension
   $\NN_0 \prec_{\sf end} \MM$. 
   The following are equivalent$:$
   
   \begin{itemize}
   
   \item[(1)] $(\MM, \cod(\NN_0 / M))$ is $\aleph_0$-strong.
   
   \item[(2)] There is a minimal extension $\NN \succ_{\sf end} \MM$ such that \\ 
   $\cod(\NN / M) = {\mathfrak X}$. \qed
   
   \end{itemize}}
   
    \bigskip

   The next lemma shows that  for many a cut $I$ of $\MM$, there is a close connection between $\dcf(I)$ and those $\kappa$ for which $(\MM|I, \cod(\MM/I))$ is  $\kappa$-strong.

    \bigskip
   
   {\sc Lemma 5.3:} {\em Suppose that $\MM \models \pa$,   $\kappa$ is an infinite cardinal and $\MM$ is $\kappa$-generated over  the  multiplicatively closed cut $I$. Then  $\dcf(I) \leq \kappa$ iff $(\MM|I, \cod(\MM/I))$ is $\kappa$-strong.}

   \bigskip
   
   {\it Proof}. Let $\MM, I$ and $\kappa$ be as given in the lemma.
   
   \smallskip
   
   $ (\Longleftarrow )$: Assume that $(\MM|I, \cod(\MM/I))$ is  $\kappa$-strong and, for a contradiction,  that  $\dcf(I) > \kappa$. Since $\MM$ is $\kappa$-generated over $I$,  there is $c \in M$ such that if 
   $C = \{(c)_i : i \in I\} \backslash I\}$, then $\inf(C) = I$. Let $A = \{\langle i,j \rangle \in I : (c)_i = j\}$. 
   Then $A \in \cod(\MM / I)$. Let $D = \{i \in I : \langle i,j \rangle \in A$ for some $j \in I\}$. Then $D \in \Sigma_1^0$-$\Def(\MM|I, \cod(\MM/I))$, and its complement $I \backslash D \in \Pi_1^0$-$\Def(\MM|I, \cod(\MM/I))$. By $\kappa$-strongness, let $\{B_\alpha : \alpha < \kappa\} \subseteq  \Sigma_1^0$-$\Def(\MM|I, \cod(\MM/I))$ whose union is $I \backslash D$. Let $C_\alpha = \{(c)_i : i \in B_\alpha\}$. 
   Thus, each $C_\alpha \subseteq M \backslash I$, and their union is $C$. Since $B_\alpha \in  \Sigma_1^0$-$\Def(\MM|I, \cod(\MM/I))$, there is $c_\alpha \in M$ such that $C_\alpha = \{(c_\alpha)_i : i \in I\}$, so that $\inf(C_\alpha) \neq I$. For each $\alpha < \kappa$, let $b_\alpha$ be such that $I < b_\alpha < C_\alpha$. Then,  
   $\inf(\{b_\alpha : \alpha < \kappa\}) = I$, contradicting that $\dcf(I) > \kappa$. 
   
   \smallskip
   
   $(\Longrightarrow)$: Assume $\dcf(I) \leq \kappa$ and that $\inf(\{b_\alpha : \alpha < \kappa\}) = I$. 
   Let $X \in \Pi_1^0$-$\Def(\MM|I, \cod(\MM/I))$. There is $F \in \cod(\MM / I)$ such that $X = \{ x \in I : \langle x,y \rangle \not\in F$ for all $y \in I\}$. We can assume that if $x\in I \backslash X$, then there is a unique $y \in I$ such that $\langle x,y \rangle \in F$. 
   Let $A \in \Def\MM)$ be such that $A \cap I = F$. We can assume that for all $x \in M$, there is at most one $y \in M$ such that $\langle x,y \rangle \in A$. Thus, we can think of $A$ as a function. 
   For each $\alpha < \kappa$, let $B_\alpha = \{ x \in \dom(A) \cap I: A(x) \geq b_\alpha\}$. Clearly, each    $B_\alpha \in \Sigma_1^0$-$\Def(\MM|I, \cod(\MM/I))$ and $X$ is the union of the $B_\alpha$'s. Thus, $(\MM | I, \cod(\MM / I))$ is $\kappa$-strong. \qed

\bigskip

By the previous lemma,  (1) of Theorem~3.3 is equivalent to: 
 $I$ is exponentially closed and $(\MM | I, \cod(\MM / I))$ is $\aleph_0$-strong. 

It is possible for a model $\MM$ to be countably generated over a cut $I$ and still have $\dcf(I)$ being very large, as the next proposition shows.

   \bigskip
   
   {\sc Proposition 5.4}: {\em Suppose that $\MM \models \pa$ and $\kappa$ is a regular infinite cardinal. 
   There is $\NN \equiv \MM$ and an elementary cut $I$ of $\NN$ such that $\NN$ is finitely generated over $I$ and $\dcf(I) = \kappa$.}
   
   \bigskip
   
   {\it Proof.} Assume that $\MM$ is countable. Let ${\mathfrak X} \subseteq {\mathcal P}(M)$ be such  that ${\mathfrak X} \supseteq \Def(\MM)$, ${\mathfrak X}$ is countable and $(\MM, {\mathfrak X}) \models \wkl+ \neg\aca$. By Corollary~2.8, let $\MM_1 \succ_{\sf end} \MM$ be a finitely generated extension such that $\cod(\MM_1 / M) = {\mathfrak X}$. 
   Then $(\MM_1 |M, \cod(\MM_1 / M))$ is $\aleph_0$-strong but $(\MM_1 |M, \cod(\MM_1 / M)) \models \neg\aca$. (If $\kappa = \aleph_0$, we can stop right now. So, assume that $\kappa$ is uncountable.)    Let $(\NN_1, I) \equiv (\MM_1,M)$ be sufficiently saturated. With care, we can arrange that 
   $(\NN_1 |I, \cod(\NN_1 / I))$ is $\kappa$-strong but not $\lambda$-strong for any $\lambda < \kappa$. By Lemma~3.4, $\dcf(I) = \kappa$. Thus, there is $a \in N_1$ such that $\inf(\{ (a)_x : x \in I\} \backslash I) = I$. Let $\NN$ be the elementary substructure of $\NN_1$ generated by $I \cup \{a\}$. Then, $\NN$ and $I$ are as required. \qed

   \bigskip 

The next theorem characterizes when the extension $\NN \succ_{\sf cf} \MM$ in Theorem 2.1 can be non-filling (except for the case that ${\mathfrak X} = \cod(\MM / I)$, which is taken care of by Theorem~3.3).

\bigskip

{\sc Lemma 5.5}: {\em Suppose that $\MM \models \pa$,  $I$ is an exponentially closed cut $I$, $ \cod(\MM / I) \subsetneq {\mathfrak X} \subseteq {\mathcal P}(I)$  and $\MM$ has a non-filling ${\mathfrak X}$-extension. Then$:$
\begin{itemize}

\item [\rm{(a)}] If $I$ is strong, then $(\MM | I, {\mathfrak X})$ is $\aleph_0$-strong \hspace{-3pt}$;$ and

\item [\rm{(b)}] if  $\Delta_2^0$-$\Def(\MM |I, \cod(\MM / I)) \subseteq {\mathfrak X}$, then $I$ is strong.

\end{itemize}}

\bigskip

{\it Proof}.  Let $\NN \succ_{\sf cf} \MM$ be a  non-filling ${\mathfrak X}$-extension.

\smallskip

 (a)  Suppose that  $I$ is a strong cut of $\MM$.  By Lemma~5.3,  $\dcf^\MM(I) = \aleph_0$. Since $I = \gcis(\MM, \NN)$ and $\NN$ is non-filling, then 
$\dcf^\NN(I) = \dcf^\MM(I)$, so by Lemma~5.3 again, $(\NN|I, \cod(\NN / I))$ is $\aleph_0$-strong. 
Since $(\MM | I,{\mathfrak X}) =  (\NN|I, \cod(\NN / I))$, then $(\MM | I, {\mathfrak X})$ is ${\aleph_0}$-strong. 

\smallskip 

(b) Suppose that $I$ is not a strong cut of $\MM$. Let $a \in M$ be such that 
$\inf(\{(a)_i : i \in I\} \backslash I) = I$. Thus, for every $b \in M$ there is $i \in I$ such that 
$I < (a)_i < b$. Since $\NN$ is non-filling, then for every $b \in N$ there is $i \in I$ such that 
$I < (a)_i < b$. Let   $B = \{ i \in I : (a)_i \in I\}$. Then $B \not\in {\mathfrak X}$ and $B \in \Sigma_1^0$-$\Def(\MM | I, \cod(\MM / I))$, implying that $\Sigma_1^0$-$\Def(\MM | I, \cod(\MM / I)) \not\subseteq {\mathfrak X}$. \qed

\bigskip

{\sc Theorem 5.6}: {\em Suppose that $\MM \models \pa$,  $\MM$ is countably generated over the exponentially closed cut $I$, $ \cod(\MM / I) \subsetneq {\mathfrak X} \subseteq {\mathcal P}(I)$  and $\MM$ has a countably generated ${\mathfrak X}$-extension. The following are equivalent$:$

\begin{itemize}

\item [(1)]  $(\rm{a})$ If $I$ is strong, then $(\MM | I, {\mathfrak X})$ is $\aleph_0$-strong \hspace{-3pt}$;$ and

\noindent $(\rm{b})$ if  $\Delta_2^0$-$\Def(\MM |I, \cod(\MM / I)) \subseteq {\mathfrak X}$, then $I$ is strong.

\item[(2)] 
 There is a  countably generated, non-filling extension $\NN \succ_{\sf cf} \MM$ such that $\gcis(\MM, \NN)$ and $\cod(\NN / I) = {\mathfrak X}$. 
 
 \item[(3)] 
 There is a  finitely generated, non-filling extension $\NN \succ_{\sf cf} \MM$ such that $\gcis(\MM, \NN)$ and $\cod(\NN / I) = {\mathfrak X}$.

 \end{itemize}}

\bigskip

{\it Proof.} Let $\MM, I$ and ${\mathfrak X}$ be as given. We will prove $(3) \Longrightarrow (2) \Longrightarrow (1) \Longrightarrow (3)$. 

\smallskip

$(3) \Longrightarrow (2)$: Trivial.

\smallskip

$(2) \Longrightarrow (1)$:  By Lemma~5.5. 
 
\smallskip

$(1) \Longrightarrow (3)$:  Suppose  that $(1)$ holds, so that both (a) and (b) are true. We will consider two cases depending on whether or not $I$ is strong. 

\smallskip

{\it $I$ is not strong}: By Theorem~2.1, let  $\NN \succ_{\sf cf} \MM$ be any finitely generated ${\mathfrak X}$-extension.
By Theorem~4.1, $\NN$ is not a filling extension, so it must be non-filling. 
\smallskip

{\it $I$ is strong}: Thus, $(\MM | I, {\mathfrak X})$ is ${\aleph_0}$-strong. We will augment the proof of $(3) \Longrightarrow (1)$ of Theorem~2.1 by adding the following requirement (S3) to the two requirements (S1) and (S2) that the sequence $\Phi_0(x) \subseteq \Phi_1(x) \subseteq \Phi_2(x) \subseteq \cdots$ of $e$-big allowable sets must have. We let $e, G$ and ${\mathfrak X}_0$ be as in that proof.

 \begin{itemize}

  \item[(S3)] For each $1$-ary Skolem ${\mathcal L}(M)$-term $t(x)$, there are $n < \omega$, $a \in I < b \in M$ 
  and $\varphi(x) \in \Phi_n(x)$ such that 
  $$
  \MM \models \forall x[\varphi(x) \into (t(x) \leq a \vee t(x) \geq b)].
  $$

\end{itemize}
It is clear that if (S3) is satisfied, then $\NN$ does not fill $I$. 

Suppose that we are at the stage at which we have just obtained  $\Phi_m(x)$ and are considering 
$t'(x,u)$, which is one of the countably many $2$-ary Skolem ${\mathcal L}(G)$-terms.  
Let $d > I$ be a bound for $\Phi_m(x)$. Our goal is to arrange that whenever $i \in I$ and $t(x) = t'(x,i)$, then the conclusion of (S3) holds. By Lemma~2.1.3, we can assume that $\Phi_m(x)$ is tree-based and that 
$\Phi_m(x) = \{\theta(x,s) : s \in B\}$, where $\theta(x,u)$ is a tree-formula and $B \in {\mathfrak X}$ is a 
$2^{<I}$-path. We further assume that $\MM \models \forall x,u[\theta(x,u) \into x < e]$. For each $s \in 2^{<M}$, let $X_s$ be the set defined by $\theta(x,s)$. 
Thus, $\MM \models |X_s| \geq d$ for each $s \in B$. Without loss, we can assume that 
$\MM \models \forall s \in 2^{<M}[ |X_s| \geq d]$. 

We can assume (by (S2)) that there is $A \in {\mathfrak X}$ such that the formula
$$
 \exists u,w[v = \langle u,w \rangle \wedge t'(x,u) \leq w]
$$ 
represents $A$ in $\Phi_m(x)$. 

Let $A$ be the set of $i \in I$ such for every $a \in I$ there is $\varphi(x)$ in $\Phi_m(x)$ such that 
$$
 \MM \models \forall x[\varphi(x) \into t(x,i) > a]. 
$$ 
Clearly, $A \in \Pi_1^0$-$\Def(\MM | I, {\mathfrak X})$, so there are countably many sets $A_0,A_1,$ $ A_2, \ldots$ in 
$\Sigma_1^0$-$\Def(\MM | I, {\mathfrak X})$ such that $A = \bigcup_{k<\omega}A_k$. (If $A = \varnothing$, then for every $i \in I$ there are $a \in I$ and $\varphi(x) \in \Phi_m(x)$ such that $\MM \models \forall x[\varphi(x) \into t'(x,i) \leq a]$, so we are done.)  
For each $k < \omega$, we will consider $A_k$ at some future stage. For notational simplicity, let's suppose that we are considering $A_k$ right now. 

For $y \in 2^{<M}$, let 
$$
C(y) = \{ j \in M : \langle i,j \rangle < \ell(y) {\mbox{ and }}(y)_{\langle i,j \rangle} = 0 {\mbox{ for some }} i \in M\}.
$$
If $z \lhd y \in 2^{<M}$, then $C(z) \subseteq C(y)$. 
Since $A_k \in \Sigma_1^0$-$\Def(\MM | I, {\mathfrak X})$, there is $Y \in {\mathfrak X} \cap 2^{<I}$ such that 
$A_k = \bigcup_{y \in Y}C(y)$. 

Let $\theta'(x,u,y)$ be the formula
$$
\theta(x,u) \wedge \forall j \in C(y) [t(x,j) \geq j],
$$ 
and let $X_{s,y}$ be the set defined by $\theta'(x,s,y)$. 
Clearly, if $s \in B$ and $y \in Y$, then $|X_{s,y}| >I$. For each $s \in B$ let $r_s \in M$ be such that 
$\MM \models |X_{s,y}| = r_s$, where $y \in Y$ is such that $\ell(y) = \ell(s)$. By (b), there is $r \in M$ such that $I < r < r_s$ for all $s \in B$. 

Let $\Phi'(x) = \Phi(x) \bigcup \{\theta'(x,u,y) \wedge \ell(s) = \ell(y) \wedge \forall j \in C(y)\big[t(x,j) \geq r\big]$. Then,
$\Phi'(x) \supseteq \Phi(x)$, $\Phi'(x)$ is bound by $r$, and for all $j \in A_k$ there is $\varphi(x) \in \Phi'(x)$ such that $\MM \models \forall x[\varphi(x) \into t(x,j) \geq r]$.   
\qed

\bigskip

The previous theorem required that $\MM$ have a countably generated ${\mathfrak X}$-extension. But having an $\aleph_1$-generated ${\mathfrak X}$-extension 
suffices for getting non-filling $\aleph_1$-generated extension. 

\bigskip

{\sc Corollary 5.7}: {\em Suppose that $\MM \models \pa$,  $\MM$ is countably generated over the exponentially closed cut $I$, $ \cod(\MM / I) \subseteq {\mathfrak X} \subseteq {\mathcal P}(I)$  and $\MM$ has an 
$\aleph_1$-generated ${\mathfrak X}$-extension. The following are equivalent$:$

\begin{itemize}

\item [(1)]  $(\rm{a})$ If $I$ is strong, then $(\MM | I, {\mathfrak X})$ is $\aleph_0$-strong \hspace{-3pt}$;$ and

\noindent $(\rm{b})$ if  $\Delta_2^0$-$\Def(\MM |I, \cod(\MM / I)) \subseteq {\mathfrak X}$, then $I$ is strong. 

\item[(2)] 
 There is an ${\aleph_1}$-generated, non-filling extension $\NN \succ_{\sf cf} \MM$ such that $\cod(\NN / I) = {\mathfrak X}$. 
 
 \end{itemize}}
 
 \bigskip
 
 {\it Proof}. Let $\MM, I$ and ${\mathfrak X}$ be as given. Then $(2) \Longrightarrow (1)$ by Lemma~5.5. 
 The converse $(1) \Longrightarrow (2)$ follows from  Theorem~3.3 in the same way that Corollary~2.10 
 follows from Theorem~2.1. \qed
 
  \bigskip
  
  The next corollary improves Corollary~3.5.
  
  \bigskip
  
  {\sc Corollary 5.8}: {\em Suppose that $\MM \models \pa$,  $\MM$ is countably generated over the exponentially closed cut $I$, $\dcf(I) = \aleph_0$,  $\cod(\MM / I) \subseteq {\mathfrak X} \subseteq {\mathcal P}(I)$  and $\MM$ has a countably generated, non-filling ${\mathfrak X}$-extension. Whenever $\kappa \geq \lambda = |I|$, there is a non-filling $\NN \succ_{\sf cf} \MM$ such that 
$I = \gcis(\MM, \NN) = I^\NN_\lambda$, $J^\NN(\lambda)= \kappa$ and $\cod(\MM/I) = \cod(\NN/I)$.} 

\bigskip

{\it Proof}. Let $\NN_0 \succ_{\sf cf} \MM$ be a finitely (or countably) generated, non-filling ${\mathfrak X}$-extension. Then apply Corollary~3.5 to get a non-filling $\NN \succ_{\sf cf} \NN_0$ such that 
$I = \gcis(\NN_0, \NN) = I^\NN_\lambda$, $J^\NN(\lambda)= \kappa$ and $\cod(\NN_0/I) = \cod(\NN/I)$.
\qed

\bibliographystyle{plain}

\end{document}